\documentclass[11pt,twoside]{article}

\pdfoutput=1

\usepackage{latexsym}
\usepackage{amssymb,amsbsy,amsmath,amsfonts,amssymb,amscd}
\usepackage{graphicx}
\usepackage{hyperref}
\usepackage{setspace}
\usepackage{pdfsync}
\usepackage{color}
\newcommand {\wt}[1] {{\widetilde #1}}
\newcommand{\commentout}[1]{}
\newcommand{\R}{\mathbb{R}}

\newcommand {\al} {\alpha}
\newcommand {\e}  {\varepsilon}

\newcommand {\sg} {\sigma}

\newcommand {\lb} {\lambda}
\newcommand {\Chi} {{\bf \raise 2pt \hbox{$\chi$}} }
\newcommand {\f}   {\frac}
\newcommand {\p}   {\partial}
\newcommand{\dis}{\displaystyle}
\newcommand {\proof} {\noindent {\bf Proof}. }
\newcommand{\re}{\eqref}
\newcommand{\beq}{\begin{equation}}
\newcommand{\eeq}{\end{equation}}
\newcommand{\bea} {\begin{array}{rl}}
\newcommand{\eea} {\end{array}}
\newcommand{\bepa}{\left\{ \begin{array}{l}}
\newcommand{\eepa} {\end{array}\right.}
\newtheorem{theorem}{Theorem}[section]
\newtheorem{lemma}[theorem]{Lemma}

\newcommand{\qed}{{ \hfill
                       {\unskip\kern 6pt\penalty 500 \raise -2pt\hbox{\vrule\vbox to 6pt{\hrule width 6pt
                       \vfill\hrule}\vrule} \par}   }}

\title{\Large \bf Populational adaptive evolution, chemotherapeutic resistance and multiple anti-cancer therapies}

\author{
Alexander Lorz\thanks{UPMC Univ Paris 06, CNRS UMR 7598, Laboratoire Jacques-Louis Lions, 4, pl. Jussieu F75252 Paris cedex 05.  Email: lorz@ann.jussieu.fr} \footnotemark[4]
\and Tommaso Lorenzi\thanks{Department of Mathematics, Politecnico di Torino, Corso Duca degli Abruzzi 24, I10129 Torino. Email: tommaso.lorenzi@polito.it}
\and Michael E. Hochberg\thanks{Institut des Sciences de l'Evolution, CNRS, UniversitŽ Montpellier 2, Place Eugene Bataillon, F34095 Montpellier and Santa Fe Institute, 1399 Hyde Park Rd, Santa Fe, New Mexico, USA. Email: mhochber@univ-montp2.fr}
\and Jean Clairambault\thanks{INRIA-Rocquencourt, EPI BANG. Email: jean.clairambault@inria.fr} \footnotemark[1]
\and Beno\^ \i t Perthame\footnotemark[1] \footnotemark[4] \thanks{and Institut Universitaire de France.  Email: benoit.perthame@upmc.fr}
}

\date{\today}

\begin{document}
\maketitle
\pagestyle{plain}
\pagenumbering{arabic}

\begin{abstract}
Resistance to chemotherapies, particularly to anticancer treatments, is an increasing medical concern. Among the many mechanisms at work in cancers, one of the most important is the selection of tumor cells expressing resistance genes or phenotypes. Motivated by the theory of mutation-selection in adaptive evolution, we propose a model based on a continuous variable that represents the expression level of a resistance gene (or genes, yielding a phenotype) influencing in healthy and tumor cells birth/death rates, effects of chemotherapies (both cytotoxic and cytostatic) and mutations. We extend previous work by demonstrating how qualitatively different actions of chemotherapeutic and cytostatic treatments may induce different levels of resistance.

The mathematical interest of our study is in the formalism of constrained Hamilton-Jacobi equations in the framework of viscosity solutions. We derive the long-term temporal dynamics of the fittest traits in the regime of small mutations. In the context of adaptive cancer management, we also analyse whether an optimal drug level is better than the maximal tolerated dose.
\end{abstract}

\bigskip

\\[4mm]

\section{Introduction}
\label{sec:intro}

Many pharmacotherapeutic processes, in anticancer,  antiviral, antimalarial, or antibiotic therapies may fail to control proliferation, because the target (virus, cell, parasite) population becomes resistant to the drug(s). This may occur either because of evolution independent of drug(s), or because the drug(s) select for
resistance. The occurrence of drug resistance is thus a major obstacle to therapy success.

Resistance is commonly seen in a diverse range of systems including agriculture pest resistance to toxic crops \cite{LMV, TG},  bacterial resistance to antibiotics \cite{magal}, mosquito resistance to antimalarial treatments \cite{bacaer} and references therein) and bacterial resistance to high copper concentrations  \cite{jerez}. \\
 
There is overwhelming evidence that populations of cancer cells may harbour a diverse array of chemo-resistant subpopulations.

Resistance mechanisms in cancer may be classified  according to at least two major effects:
\\
$\bullet$ {\em Cell adaptation.}  Cells may adapt metabolically  in response to drugs. A typical example is the release of resistance proteins such as P-gp and other ABC transporters (ATP Binding  Cassette), and studies on large human cohorts have shown their correlation with resistance in Acute Myeloid Leukemia  \cite{marzac}), since they promote the efflux  of  toxic molecules from leukaemic cells  \cite{jc, GFB}. Induction of ABC transporter gene transcription by stressful conditions (including exposure to anticancer drugs) is a likely contributing mechanism, as reviewed in \cite{scotto}. This resistance mechanism can further be enhanced by protein exchange  between cells \cite{rafii}, and its properties have been analysed using a mathematical model \cite{magal2}. In as much as these mechanisms occur in individual cells  and do not involve differential selection on variants resulting from mutations, we distinguish such {\em cell adaptation} mechanisms from other mechanisms stemming from {\em mutations}.
\\
$\bullet$  {\em Mutations.} It is generally accepted that tumor cells usually have higher division and higher mutation rates than normal cells; possible mechanisms for this have been investigated in bacterial populations in stressful conditions, inducing the formation of error-prone polymerase \cite{kivisaar}, from which one may infer similar explanations for cancer cells. Moreover, when a drug targets a specific DNA site or a specific protein, a mutation that deletes this site or induces a  change in protein conformation may confer resistance to mutant cells \cite{shah04}.\\
\\
To elicit the relative contributions of  these two mechanisms presents a considerable challenge.  The answer may depend on many parameters such as the type of tumor and therapy. In the latter case, such mutations (possibly by a single amino acid of the target protein) may induce resistance to the drug at stake, whereas in the former case of multiply targeted drugs (or drugs with unidentiÞed speciÞc targets) a re-wiring of intracellular signaling networks is a plausible mechanism, as discussed in \cite{McCormick}. As a first approach to this complex question, we focus in this paper on continuous mechanisms that involve resistance induced by the drug environment on cell populations, and not on mechanisms that can be represented by a constant probability of mutation resulting in drug resistance, as in classical studies such as \cite{goldie-coldman} or more recent work based on
 evolutionary models in a stochastic setting \cite{Komarova}.

In this paper we are particularly interested in the latter mechanism, which is closely related to the process of Darwinian evolution, according to which certain rare mutants may have positive growth rates and be selected in environments that would otherwise result in ecological extinction. The term of `evolutionary rescue' has been coined for situations where a population is saved from extinction by natural selection of genetic variants either existing in a population when environmental conditions deteriorate, or emerging during the otherwise inexorable decline to extinction. In this view, because some mutant cells have a fitness advantage, they will come to dominate dynamics and the population as a whole. We propose a mathematical model that takes into account birth, death and mutation rates of healthy and cancer cells depending on  the level of  resistant phenotype expression. Different expression levels would be reflected by the cell content of proteins (e.g., ABC transporters) responsible for expelling toxic molecules.  We assume this expression level to be a continuous variable, and as such amenable to provide a basis for structured population models designed according to  a  selection/mutation type of formalism previously developed using Partial Differential Equations \cite{CC0, CC, CFM2006, djmp, LMP, MW, barlesP2}.
\\

The phenomenon of drug resistance is usually described mathematically  by considering two distinct types of cells, sensitive and resistant. A recent probabilistic approach \cite{FooM} incorporated dose modulation to study the evolution of resistance under various dosing strategies. Unlike probabilistic models, we do not consider that a cell or a cell population is necessarily either totally sensitive or totally resistant to a given drug; rather we introduce a continuous physiological variable describing resistance between $0$ (completely sensitive population) and $+\infty$ (completely resistant population). Such a variable, which may be genetic or epigenetic \cite{TroyDay}, indeed provides the basis for a physiologically structured model of resistance.

Mathematical models  \cite{GSGF} incorporating phenotypic and microenvironmental factors suggest that adaptive treatment strategies can improve long-term survival compared to a maximum dose strategy. This may be achieved at the cost of maintaining the persistence of sensitive tumor cells, which are fitter than the resistant cells in low drug pressure conditions (i.e., not trying to eradicate them, an apparently paradoxical view (`adaptive therapy') \cite{Gatenby, GSGF}, most likely applicable to slowly developing tumors only).
Conversely, in cases of sudden rapid proliferation (e.g., acute leukaemias), only the employment of maximum tolerable doses is thought to have chances of controlling or eradicating the cancer cell population.

Furthermore, depending on the nature of the tumor, different sorts of drugs can be used either separately or in combinations, in order to reduce the probability of resistance emergence and reduce side-effects on healthy cells \cite{siga}. This aspect is studied in \cite{levy1, levy2} using ordinary differential systems which give a macroscopic understanding on the use of multitherapies. In particular, besides cytotoxic drugs, here we focus on an additional class of therapeutic agents,  so-called cytostatic drugs, which act by slowing down cancer cell proliferation and tumor growth. Cytostatic drugs have lower toxicity for healthy cells and reduce the emergence of resistance, which usually follows from treatments with cytotoxic drugs. In fact, they allow the survival of a small number of chemosensitive cells, which can reduce the growth of resistant clones through competition for space and resources. Moreover, whereas cytotoxic drugs remove sensitive cells from the population, they do not prevent the proliferation of chemo-resistant cells. There are two important effects of such therapies on cancer resurgence and its aggressiveness. Firstly, proliferating resistant cells may produce a series of increasingly aggressive subpopulations through natural selection \cite{bozic}. Secondly, toxic effects of the drug on chemoresistant cells, although not lethal, may induce genetic instabilities resulting in the generation of more aggressive mutant cancer cell lines \cite{friedberg}.  These concerns within an expanding chemoresistant population are reduced under combined chemostatic and chemotoxic therapies, since toxicity of the latter is reduced, and the chemostatic drug holds in check the multiplication of the chemoresistant subpopulation.
\\

The paper is organized as follows. We first present a model for healthy and tumor cells, structured by a variable representing the expression level of the resistance gene. We consider cases with and without chemotherapy. We also analyze the global behavior of solutions without mutations. Then, in Section \ref{sec:aasm}, we turn to the analysis of situations with small mutations and employ the methodology of populational adaptive evolution to describe the dynamics of the fittest cells. Section \ref{sec:opt} is devoted to the question of analyzing, for this simple model, whether the maximal tolerated dose is the optimal choice to impede the evolution of resistance. Section \ref{sec:numerics} presents several numerical simulations that illustrate our theory. Finally, we propose in Section \ref{sec:model2} a model for the dynamics of healthy and cancer cells exposed to cytotoxic and cytostatic drugs. The model relies on the same structured population formalism presented at the beginning of the paper and it is analyzed by means of numerical simulations looking for combined doses resulting in treatment optimization.

\section{A structured population model for healthy and tumor cells under the effects of cytotoxic drugs}
\label{sec:model}

A major constraint in chemotherapy is to keep the toxic effect on healthy cells below a critical threshold.  For that reason, it is preferable for a model to be able to represent both healthy and tumor cells. This is what we propose here in a simple equation that can take into account both types of cells, thus tackling the two major issues of medical treatments in general: adverse toxic effects in healthy cells and occurrence of drug resistance in diseased cells. At this stage we do not represent the spatial extension of tissues; this is an ingredient that can be included in a second stage. It is indeed more important from a therapeutic point of view to structure the two populations according to a physiological variable, e.g., a continuously evolving genetic trait or a cell content in specific proteins, that physiologically describes the evolution of the population as influenced by the treatment, than to include a space variable in the model.

\subsection{Selection/mutation model}
\label{sec:selection}

Let $n_H(x,t)$ / $n_C(x,t)$ denote the population density of healthy / cancer cells with gene resistance expression level $x$ at time $t$. In the sequel we will use the term `gene expression' meaning not only expression of one supposed resistance gene, but more generally of several genes yielding together a continuous drug resistance phenotype, more or less in the same way as others represent evolution towards malignancy in colorectal cancer \cite{FearonVogelstein}, revisited by  \cite{Sprouffske}. Note that in the model we develop below, we do not represent evolution towards malignancy, assuming the existence of an already constituted cancer cell population (with subscript $C$), opposed to a healthy cell population (with subscript $H$). The resistance level $x$ to the drug can be measured either by the average molecular cell concentration or activity of ABC transporters that are known to be associated with resistance to the drug, or by the minimal drug concentration to kill the population (more precisely the LD50, LD90... of pharmacologists, i.e., the minimal dose to kill a given percentage of the cell population). Different strains of the same initial cell lineage, selected by exposure to progressively increasing doses of the same drug, are indeed available in tumor banks, such as acute myeloblastic leukemia (AML) cells at the Tumor Bank of St. Antoine Hospital in Paris \cite{Zhou}. We {\em a priori} take $0 \leq x < \infty $, even though the design of the model should lead to limitations on the possible values of $x$.
\\

The growth dynamics of healthy and tumor cells with a chemotherapy is given by the system
\beq\begin{array}{rl}
\dis\f{\p}{\p t} n_H(x,t) =  &\Big [  \overbrace{\f{1-\theta_H}{\big(1+\rho(t) \big)^\beta} \; r(x)}^{\rm growth \; with \; homeostasis}  -  \overbrace{d(x)}^{\rm natural \; apoptosis} -  \overbrace{c(t) \mu_H(x) }^{\rm effect \; of \; drug} \Big]  n_H(x,t)
\\[5mm]
& + \dis\f{\theta_H}{ \big(1+\rho(t)\big)^\beta} \; \underbrace{ \int r(y) M_{\sg_H}(y, x) n_H(y,t) dy}_{\rm birth\; with \; mutation},
\end{array}
\label{eq:healthy}
\eeq
\beq\begin{array}{rl}
\dis\f{\p}{\p t} n_C(x,t) =  &\Big [   (1-\theta_C) \; r(x)  -  d(x) - c(t) \mu_C(x) \Big]  n_C(x,t)
\\[5mm]
& +\theta_C \; \dis \int r(y) M_{\sg_C}(y, x) n_C(y,t) dy,
\end{array}
\label{eq:cancer}
\eeq
and the total population is defined as
\beq
\rho(t) = \rho_H(t) + \rho_C(t)  , \qquad \rho_H(t) = \int_{x=0}^\infty  n_H(x,t) dx, \quad \rho_C(t) = \int_{x=0}^\infty  n_C(x,t) dx.
\label{eq:rho}
\eeq

The following notations, interpretations and assumptions are used:
\\[2mm]
$\bullet$ \;   $r(x)$ denotes the basic  reproduction rate and $d(x)>0$ denotes the basic death rate, which depend on the gene expression level $x$. In order to incorporate a cost to produce the resistance gene, we assume $r$ is  decreasing, $d$ is increasing
\beq
 r(0 )> d(0) >0, \qquad  r'(Ê\cdot) <0, \qquad r(+\infty) =0, \qquad  d'(Ê\cdot)  > 0.
\label{as1}
\eeq
$\bullet$ \; $0 \leq \theta_{H,C} <1$ denotes the proportion of divisions with mutations and we can assume it to be higher for cancer cells.
\\
$\bullet$ \; $\beta >0$ is introduced, with the simplest possible form, to impose healthy tissue homeostasis (see below). For tumor cells, uncontrolled proliferation is obtained by birth terms that do not depend on the total population (i.e., we do not consider density-dependent inhibition in cancer cells).
\\
$\bullet$ \; $c(t)$ denotes the dose of chemotherapy. Here we assume it has only an effect on increasing apoptosis.
\\
$\bullet$ \; $\mu_{H,C}(x) $ represents the  phenotypically dependent response to the drug; drugs are designed to target cancer cells more than healthy cells. The effect of drugs is here assumed to be summed up directly on mortality (i.e., in this simple setting, not involving the cell division cycle, we do not consider drug effects at cell cycle phase transitions in proliferating cell populations, see e.g., \cite{KS} for a discussion on this point) with a rate $\mu$ that depends on the gene expression $x$ (sensitivity to the drug) and we assume the properties
\beq
\mu_{H,C}(\cdot)>0, \qquad \f{d \mu_{H,C}}{dx}  <0, \qquad \mu_{H,C}(\infty) = 0, \qquad \mu_{H}< \mu_{C} .
\label{as2}
\eeq
This means that the drug always decreases the growth rate but its efficiency is lower on cells expressing a higher resistance gene expression.
Also because the drug is efficient against cancer cells we assume that $c$ is such that
\beq
r(0)-d(0) - c\mu_{C}(0) < 0 .
\label{as3}
\eeq
\\
$\bullet$ \; $M_{\sg_{H,C}}(y,x) \geq 0$ denotes the probability that a mutation of a cell with gene expression $y$ leads to a daughter cell with level $x$ and $\sg_{H,C}>0$ measures the average size of these mutations. This means that
$$
 \int_{x=0}^\infty x M_\sg(y,x) dx=\sg.
$$
We can assume that this size is larger in tumor than in healthy cells, in which control of the genome is more strictly ensured and mutations, if they occur, are of little consequence; hence $\sg_C> \sg_H$. To be consistent with this interpretation, we impose (further conditions are imposed later on)
\beq
 \int_{x=0}^\infty  M_\sg(y,x) dx=1 .
\label{eq:asmut}
\eeq

It is useful to write the equations on $n_C$ and $n_H$ as
$$
\f{\p}{\p t} n(x,t) = n(x,t) R(x, t) +\theta \;  \int r(y) M_{\sg}(y, x) [n(y,t)-n(x,t)] dy,
$$
so as to better see that the {\em net growth rate} (fitness) of the cells is given by
\beq \begin{cases}
R_H(x, t) & = \dis\f{1}{\big(1+\rho(t) \big)^\beta} \;  r(x) - d(x) - c \mu_H(x),
\\[4mm]
R_C(x, t) & =  r(x) - d(x) - c \mu_C(x).
\end{cases}
\label{eq:ngr}
\eeq

\subsection{Healthy or tumor cells without mutations}
\label{sec:htc}

In order to motivate the introduction of these different terms, we will prove in the next sections  that the model solutions exhibit behaviors that are in accordance with the biomedical interpretation. It is simpler to state results when there are no mutations, which we do now and postpone to the next sections the case with mutations.
\\

We introduce our assumptions and notations in each one of these three cases of interest independently: healthy cells only, cancer cells without therapy, resistance in cancer cells generated by therapy.
\\
\\
{\bf \large Healthy tissue, no therapy, no tumor cells.} In a healthy tissue we take $n_C \equiv 0$. Then, in the absence of therapy (i.e., $c(t)=0$), the classical theory, that we recall below, shows that as $\sg$ vanishes, the population model always selects the gene expression $x=0$. In other words, $x=0$ is a globally attracting Evolutionary Stable Distribution:
\begin{lemma} [Healthy cells only, no mutations] We assume \eqref{as1} and that all data are Lipschitz continuous; we also assume
$$
\theta_H=0, \qquad n_H^0 >0, \qquad n^0_H \in L^1, \qquad n_C \equiv 0.
$$
Then, we have
$$
n_H(x,t) \underset{t \to \infty}{\longrightarrow} \bar \rho_{H,\infty} \delta (x),
$$
weakly in the sense of measures where the number density of cells $ \bar \rho_{H,\infty}$ defined by zero growth rate (homeostasis) is computed as
\beq
\dis \f{1 }{\big(1+\bar \rho_{H,\infty} \big)^\beta} \; r(0)  =   d(0) .
\label{aa:rhoinfty}
\eeq
\label{lm:h}
\end{lemma}
In words, with no mutations to create variability, the monoclonal  population is globally attractive with the fittest trait $x=0$.

The proof of this result follows from a general analysis that can be found in \cite{LMP, Pe, barlesP2}. It uses an estimate on the total variation of $\rho_H(t)$ that holds in a very general framework that includes the case at hand.
\\

Keeping mutations, a similar conclusion can be reached but in the framework of asymptotic analysis; to show this is a job that uses more elaborate ingredients and that is explained in Section \ref{sec:aah}.
\\
\\
{\bf \large Cancer cells, no  therapy.}
For healthy cells, the total cell number $\rho(t)$ induces a limitation on the multiplication of cells that represents homeostasis (i.e., the capacity of a tissue to stop growing, then differentiate and commit itself to fulfill a given physiological task). For cancer cells, we observe an uncontrolled growth in the absence of therapy whatever the gene expression level (at least close enough to $x\approx 0$, otherwise the cost to produce it can compensate the birth/death ratio).  However, they remain with the lowest resistance gene because of the reproduction advantage when $x = 0$. This can be derived from the model since we have the
\begin{lemma}  [Cancer cells only, no mutations] When $c(t)\equiv 0$ and $\theta_C$ vanishes, solutions to \eqref{eq:cancer} with an initial data $n_C^0(x) >0$ near $x\approx 0$ satisfy, with an exponential rate
$$
\rho_C(t) \underset{t \to \infty}{\longrightarrow} \infty, \mbox{ and weakly in the sense of measures } \f{n_C(x,t)}{ \rho_{C}(t)} \underset{t \to \infty}{\longrightarrow} \delta (x) .
$$
\label{lm:tumor}
\end{lemma}

\proof From assumption \eqref{as1}, there is a value $x_c>0$ for which
$$
r(x) -d(x) \geq g_m := \f 12 [r(0) -d(0)]  >0, \quad \forall x \in [0,x_c].
$$

We write after integration in $x$
$$
\f{d}{dt} \int_0^{x_c} n_C(x,t) dx = \int_0^{x_c} [r(x) -d(x)] n_C(x,t) dx \geq g_m  \int_0^{x_c} n_C(x,t) dx.
$$
Therefore the first statement is proved since
$$
 \int_0^{x_c} n_C(x,t) dx \geq  e^{\lb_m t} \int_0^{x_c} n_C^0(x) dx, \qquad \lb_m =g_m.
$$

We also notice that the same argument, using a small interval closer to  $x\approx 0$ tells us that  for all $\lb < r(0) -d(0)$, there is a $x_\lb$
such that
$$
 \int_0^{x_\lb} n_C(x,t) dx \geq  e^{\lb  t} \int_0^{x_\lb} n_C^0(x) dx,
$$
and thus, for all $\lb < r(0) -d(0)$, there is a value $\rho(\lb)$ such that
$$
\rho_{C}(t) \geq \rho(\lb) e^{\lb t}.
$$
From this, we easily conclude the second statement which is a consequence of Laplace formula. Indeed, for a given $y$, we can choose from assumption  \eqref{as1}, $\lb > r(y)- d(y) $ and thus for $x\geq y $
$$
\f{n_C(x,t)}{\rho_{C}(t)} \leq  n_C^0(x)  \f{e^{[r(x)-d(x)]t} }{\rho(\lb) e^{\lb t} } \underset{t \to \infty}{\longrightarrow} 0.
$$
In other words, the cell density is asymptotically small for all $x>0$, and thus the solution concentrates at $x=0$.
\qed

The case with mutations is treated in Section \ref{sec:aac}.

\medskip
\noindent {\bf \large Cancer cells with therapy} This is the most interesting situation where resistance may follow from selection of cells with a high level of gene expression. This can arise because the drug concentration $c$ is limited in order to keep healthy cells below a given toxicity threshold.
 This situation occurs if the maximum net growth rate  of tumor cells, as defined through \eqref{eq:ngr},  is positive for some $x>0$. Here we concentrate on the simplest case of a constant therapy and we make the additional assumption that, for all level of $c$, there is a maximal resistance gene expression $x_c>0$ which maximizes the fitness
\beq
 \text{for } x \neq x_c, \quad r(x)- d(x) - c \mu_C(x) <  \overline{R_c} := r(x_c)- d(x_c) - c \mu_C(x_c) >0.
\label{as:xr}
\eeq
Notice however that
$$
 \overline{R_c} \leq r(x_c)- d(x_c) \leq r(0) -d(0) = R_0.
$$

With this assumption, we can establish the
\begin{lemma}  [Cancer cells with therapy, no mutations] Assume $c \equiv 1$, $\theta_C=0$ and \eqref{as:xr}. Then, the solutions to \eqref{eq:cancer} with an initial data $n_C^0(x) >0$ satisfy, with an exponential rate
$$
\rho_C(t) \underset{t \to \infty}{\longrightarrow} \infty, \qquad \f{n_C(x,t)}{ \rho_{C}(t)} \underset{t \to \infty}{\longrightarrow} \delta (x-x_c) .
$$
\label{lm:therapy}
\end{lemma}

We do not repeat the proof of this result which is the same as for Lemma \ref{lm:tumor} and gives that $\rho_C(t) e^{\lb t}\underset{ t \to \infty }{\longrightarrow} \infty$ for all $\lb < R_c$.

\section{Asymptotic analysis with small mutations}
\label{sec:aasm}

A clear mathematical way to precisely express the results mentioned above, i.e., that a specific level of gene expression is selected, is through asymptotic analysis for small mutations and observing the dynamics in the long run. In ecology, these models generate populations that are highly concentrated on some well separated traits with possible branching and polymorphism. For the case at hand, the mutation rates might be too high to reduce the analysis to this limit, but it is still a way to express mathematically what happens for small, but non zero, frequencies of mutation.

In this section we state the behaviors that can be derived from analysis in the different cases already mentioned without mutations in Section  \ref{sec:htc}. They give details on the evolution of the population towards the ESD (Evolutionary Stable Distribution, see \cite{JR}), which is that at each time a distribution close to a Dirac mass is obtained.

Even though these are not the most interesting cases, the results listed below are also valid when there are no mutations ($\theta_H=0$ or $\theta_C=0$) and this is the case used later in the numerical simulations of the solutions.

\subsection{Asymptotic analysis: case of healthy cells}
\label{sec:aah}
With the aim of studying the dynamics of the model in the aforementioned limit of small mutations and many generations, we introduce a small parameter $\e$, which is used to rescale time and define a more specific choice of the mutation kernel, so that the population dynamics  \eqref{eq:healthy} for healthy cells only and without treatment can be rewritten in the form
\beq\left\{\begin{array}{rl}
\e \dis\f{\p}{\p t} n_H(x,t) =  &\Big [ \dis \f{1-\theta_H}{\big(1+\rho_H(t) \big)^\beta} \; r(x)  -   d(x) \Big]  n_H(x,t)
 \\[4mm]
 & \qquad  \qquad   \qquad + \dis\f{ \theta_H}{ \big(1+\rho_H(t)\big)^\beta} \;  \int r(y)  \f 1 \e \wt M(y, \f{y-x}{\e}) n_H(y,t) dy ,
 \\[6mm]
 \rho_H(t) = & \dis \int_{x=0}^\infty n_H(x,t)  dx.
\end{array}\right.
\label{aa:healthy}
\eeq
We have written the mutation kernel in such a way that
\beq \begin{cases}
M_{\sg_H}(y, x) = \f 1 \e  \wt M(y, \f{y-x}{\e}) \qquad \text{ with }Ê\; \dis \sup_{x\geq 0} \int e^{|z|^2} \wt M(x,z) dz < \infty,
\\[4mm]
 \int_{-\infty}^{x/\e}  \wt M(x,z) dz =1 , \qquad  \int_{-\infty}^{x/\e} z \wt M(x,z) dz =0,
\end{cases}
\label{aa:mutkernel}
\eeq
so that it appears directly that $x-y= O(\e)$, in other words, mutations have a small effect in terms of the change in the $x, \; y$ variables. Notice that, in these expressions, we always have $\e z\leq x$. Therefore, the last assumption implies that cells in the state $y=0$ do not undergo mutations (they could only be to $z<0$ thus contradicting the zero integral condition). This is useful to state clear and simple results, otherwise the evolutionary drift creates  an ESS (Evolutionary Stable Strategy) (see \cite{diekmann} for an introduction to the dynamical system theory of adaptive dynamics) for a slightly positive $x$ and this can be analyzed mathematically with our method but we will not do it here for simplicity of the presentation.
\\

The analysis in \cite{BMP09,barlesP1,barlesP2}  can be used and proves that, thanks to the assumptions \eqref{as1}, the fittest population is obtained when the gene is not expressed at all, in mathematical words,  $n_H(x,t) = \bar \rho_{H,\infty} \delta (x)$ is attractive for this system.
\\

Before stating a precise result, we give assumptions. Firstly, we assume that the initial population is concentrated as a `sharp Gaussian' around a state corresponding to a non-vanishing gene expression~$\bar x^0$
\beq \begin{cases}
n_{H,\e}(x,0) =e^{u_{H,\e}^0(x)/\e}, \quad  u_{H,\e}^0(x) \underset{\e \to 0}{\longrightarrow} u_H^0(x) \; \text {locally uniformly},
\\[2mm]
\max u_H^0(x)=0=u_H^0( \bar x^0), \qquad n_{H,\e}(x,0) \underset{\e \to 0}{\longrightarrow} \bar \rho_H^0 \delta(x-\bar x^0).
\end{cases}
\label{aa:as1}
\eeq
Notice that,  as $\e$ vanishes, only mutations can deviate such a  population determined by $\bar \rho_H^0 \delta(x-\bar x^0)$  from this initial state which is  asymptotically a steady state for $\theta_H=0$  as soon as the total population satisfies
\beq
\dis \f{1 }{\big(1+\bar \rho_H^0 \big)^\beta} \; r(\bar x^0)  =   d(\bar x^0) .
\label{aa:as2}
\eeq

More precisely, with these new notations we may write
\begin{theorem}[Healthy cells only] We assume \eqref{as1}, \eqref{aa:mutkernel}--\eqref{aa:as2} and that all data are Lipschitz continuous. Then, the population is monoclonal in the limit of small mutations:
$$
n_{H,\e}(x,t) \underset{\e \to 0}{\longrightarrow} \bar \rho_H(t) \delta (x- \bar x (t)),
$$
and it returns to the fittest trait $x=0$: using the notation \eqref{aa:rhoinfty}, we have
$$
\bar x(t) \underset{t \to \infty}{\longrightarrow}  0, \qquad \bar \rho_H(t) \underset{t \to \infty}{\longrightarrow}   \bar \rho_{H,\infty}.
$$
\label{aa:th1}
\end{theorem}

\proof {\em The limit $\e \to 0$.} The proof is based on the use of the Hopf-Cole transform
$$
u_\e (x,t) = \e \ln(n_H(x,t)),\qquad  \dis\f{\p u_\e }{\p t}  = \f{\e}{n_H} \dis\f{\p n_H}{\p t} .
$$
Inserting it in  \eqref{aa:healthy}, we find the equation
$$
\dis\f{\p}{\p t} u_\e =  \dis \f{1-\theta_H}{\big(1+\rho_H(t) \big)^\beta} \; r(x)  -   d(x) + \dis\f{\theta_H}{ \big(1+\rho_H(t)\big)^\beta}
 \;  \int r(y) \wt M(y, \f{y-x}{\e}) e^{(u(y,t)-u(x,t))/\e}  dy
$$
and changing the variable $y \mapsto z$ with $y= x + \e z$, we find
$$
\dis\f{\p}{\p t} u_\e =  \dis \f{1-\theta_H}{\big(1+\rho_{H,\e}(t) \big)^\beta} \; r(x)  -   d(x) + \dis\f{\theta_H}{ \big(1+\rho_{H,\e}(t)\big)^\beta}
 \;  \int r(x + \e z) \wt M(x + \e z, z) e^{(u_\e(x + \e z,t)-u_\e(x,t))/\e}  dz.
$$

At this stage, it is useful to introduce the limit in the integral term which is given by the Hamiltonian
\beq
H(x,p)  = r(x) \int \wt M(x, z) e^{p \cdot z}  dz.
\label{eq:hamiltonien}
\eeq
Notice that from \eqref{aa:mutkernel}, we derive  the properties
$$
H(x,0)= r(x), \qquad H_{p}(x,0) = 0, \qquad H_{pp}(x,p) \geq 0.
$$
By convexity, these imply in particular
\beq
H(x,p) \geq r(x).
\label{hamiltonianpos}
\eeq

One can establish two types of uniform estimates: (i) the total number of cells $\rho_{H,\e}(t)$ is uniformly bounded in $BV(0,T)$ for all $T>0$ and  $\left(\f{d}{dt}\rho_{H,\e}(t)\right)_-$ is bounded, therefore
\beq
\rho_{H,\e}(t) \underset{\e \to 0}{\longrightarrow} \bar \rho_H(t) , \qquad   \f{d}{dt}\rho_{H}(t)\geq -C.
\label{eq:rholim}
\eeq
(ii) Uniform Lipschitz estimates can be established for this kind of models (see \cite{BMP09} for the specific case of integral mutations).
\\

Therefore, we may now pass to the limit $\e \to 0$ in the equation on $u_\e$ in the viscosity sense (see references below), and we find
\beq\bepa
\dis\f{\p}{\p t} u =  \dis \f{1-\theta_H}{\big(1+\rho_H(t) \big)^\beta} \; r(x)  -   d(x) + \dis\f{\theta_H}{ \big(1+\rho_H(t)\big)^\beta}
 \;  H(x, \nabla u(x,t)) ,
\\ [5mm]
\dis \max_{0\leq x < \infty} u(x,t) =0 \qquad \forall t >0.
\eepa
\label{aah:chj}
\eeq
This equation is the fundamental new object that describes the dynamics for $\e\to 0$. It is a constrained Hamilton-Jacobi equation, which is rather different from the usual Hamilton-Jacobi equation because the $L^\infty$ contraction property is lost for instance. Nevertheless, the correct notion of solutions are still the viscosity solutions of Crandall-Lions (see \cite{barlesbook,CIL,FS}).
\\

The solution to equation \eqref{aah:chj} is a pair $(u, \rho)$ where $\rho(t)$ is the Lagrange multiplier associated with the constraints $\max_{x} u(x,t) =0 $. There are two possible approaches to rigorously prove this limit. In the approach developed in \cite{LMP}, we consider strong assumptions which allow for uniform smoothness of $u_\e$, and then the theory can be carried much further including the long time dynamics of $\bar x (t)$. In this approach concavity assumptions are needed which are difficult to meet in the case at hand (it is an interesting problem to develop the theory).

The other approach, that we follow here, has been developed in \cite{BMP09,barlesP1,djmp,barlesP2}. We consider a weak theory with our general assumptions and viscosity solutions to \eqref{aah:chj}. The function $u$ is merely Lipschitz continuous and semi-concave \cite{barlesP2}, the best possible regularity is $\rho(t) \in BV(\R^+)$ because jumps on $\rho(t)$ and $\bar x(t)$ are possible.  This analysis for weak solutions (only possible in one dimension) gives also that, along with the dynamics, the fittest trait is characterized by
\beq
0=u(\bar x(t), t),
\qquad \qquad
0= \dis \f{1}{\big(1+\rho_H(t) \big)^\beta} \; r( \bar x(t))  -   d( \bar x(t)).
\label{aah:dyn}
\eeq
This second equality is just a consequence of $\f{\p u}{\p t}= \nabla u=0$ at a maximum point (viscosity solutions are not enough and semi-concavity is necessary here).
\\

\noindent {\em The limit $t \to \infty$.} Additionally, we  also know from \eqref{eq:rholim} that $\rho_H(t)$ is non-decreasing, therefore taking into account \eqref{aah:dyn} and the definition \eqref{aa:rhoinfty}, we find
$$
\rho_H(t) \underset{t \to \infty}{\nearrow}  \bar \rho_H \leq \bar \rho_{H,\infty}, \qquad \bar x (t) \underset{t \to \infty}{\searrow}  \bar x_\infty.
$$
It remains to identify these limits.

By contradiction, assume that $\bar \rho_H <  \bar \rho_{H,\infty}$, then  we use that
$$\bea
\dis\f{\p}{\p t} u( x, t) &=  \dis \f{1-\theta_H}{\big(1+\rho_H(t) \big)^\beta} \; r(x)  -   d(x) + \dis\f{\theta_H}{ \big(1+\rho_H(t)\big)^\beta}
 \;  H(x, \nabla u(x,t))
 \\[3mm]
 & \geq  \dis \f{1}{\big(1+\rho_H(t) \big)^\beta} \; r(x)  -   d(x)
\eea
$$
using \eqref{hamiltonianpos} for the inequality. With our assumption and the definition of $ \bar \rho_{H,\infty}$, there are intervals $[x_1,x_2]$ where the right hand side is positive for all $t$  and thus $\f{\p}{\p t} u( x, t) >0$ (uniformly) on this interval which contradicts the condition $\max_x u(x,t)=0$. This proves that $\bar \rho_H =  \bar \rho_{H,\infty}$ and concludes the proof of the theorem.
\qed

\subsection{Asymptotic analysis: cancer cells, no therapy}
\label{sec:aac}

Following the lines drawn in Section \ref{sec:aah}, we can arrive at the same conclusion that, considering cancer cells only, the gene expression $x=0$ is selected even when mutations are present but small. Again a clean mathematical path towards this direction is to rescale the model for cancer cells and, following the lines of Section \ref{sec:aah}, to rewrite it  as
\beq\begin{array}{rl}
\e \dis\f{\p}{\p t} n_C(x,t) =  &\Big [   (1-\theta_C) \; r(x)  -  d(x) \Big]  n_C(x,t)
 +\theta_C \; \dis \int r(y) \f 1\e \wt M(y, \f{y-x}{\e}) n_C(y,t) dy .
\end{array}
\label{aa:cancer}
\eeq

Because it is linear and expresses unlimited growth, this model is very different from those in Section \ref{sec:aah} for healthy cells (or in ecology). In these cases there is a   limitation by nutrients or space which controls nonlinearly the growth and the dynamics through an integral variable $\rho(t)$. Nevertheless we can use the same method to analyze again the main effect, which is that $x=0$ is a globally attractive ESD; the proof highlights how a new control unknown shows up naturally and gives the same constraint and Lagrange multiplier  in the limit as in \eqref{aah:chj}. To state the result we need again to make more precise that the initial data is supposed monoclonal at a value $\bar x_0 >0$
\beq \begin{cases}
n_{C,\e}(x,0) =e^{u_{C,\e}^0(x)/\e}, \quad  u_{C,\e}^0(x) \underset{\e \to 0}{\longrightarrow} u_C^0(x) \; \text {locally uniformly},
\\[2mm]
\max u_C^0(x)=0=u_C^0( \bar x^0), \qquad n_{C,\e}(x,0) \underset{\e \to 0}{\longrightarrow} \bar \rho^0 \delta(x-\bar x^0).
\end{cases}
\label{aa:as3}
\eeq

\begin{theorem}[Cancer cells, no therapy] We assume \eqref{as1}, \eqref{aa:mutkernel}, \eqref{aa:as3} and that all data are Lipschitz continuous. Then, the population is monoclonal in the limit of small mutations:
$$
\f{n_C(x,t)}{\int_0^\infty n_C(x,t) dx} \underset{\e \to 0}{\longrightarrow}  \delta (x- \bar x (t)),
$$
and it returns to the fittest trait $x=0$
$$
\bar x(t) \underset{t \to \infty}{\longrightarrow} 0.
$$
\label{aa:th2}
\end{theorem}
The model again is not very relevant for direct interpretation because, as we see it below, the cell number density can decrease exponentially and then it explodes exponentially. A behavior reflecting only the fact that cells cannot be naturally endowed with a gene expression $x \neq 0$; this may occur if a treatment is stopped after resistance has been triggered (see next section).
\\

\proof We define the probability density
$$
p_\e(x,t) := \f{n_C(x,t)}{\rho_C(t)}, \qquad  \rho_C(t)  = \int_0^\infty n_C(x,t) dx.
$$
Because $ \rho_C(t)  $ satisfies the equation
$$
\e \f{d  \rho_C(t) }{dt} =  \int_0^\infty [r(x) -d(x)] n_C(x,t) dx
$$
we also find a closed equation for $p_\e$, namely
$$
\e \dis\f{\p}{\p t} p_\e(x,t) = p_\e(x,t) [(1 - \theta_C) r(x)-d(x)] - p_\e(x,t) I_\e(t)
+  \theta_C \int_0^\infty r(y) p_\e(y) \f 1\e \wt M(y, \f{y-x}{\e}) dy
$$
with
\begin{equation}\label{Ie}
I_\e(t) :=  \int_0^\infty p_\e(y,t) (r(y)-d(y)) dy.
\end{equation}
As mentioned earlier, this integral comes naturally to play the role of the population density $\rho(t)$ in the case of healthy cells as a Lagrange multiplier and we are back to the general setting in  \cite{BMP09,barlesP1, LMP,barlesP2}.

We can perform the Hopf-Cole transform again $u_\e = \e \ln p_\e$ and write
$$
\dis\f{\p}{\p t} u_\e = (1-\theta_H) r(x)  -   d(x) - I_\e(t) + \theta_H \;  \int r(x + \e z) \wt M(x + \e z, z) e^{(u_\e(x + \e z,t)-u_\e(x,t))/\e}  dz.
$$
We can finally pass to the limit and find, still with the Hamiltonian defined in \eqref{eq:hamiltonien}, the constrained Hamilton-Jacobi equation
\beq\bepa
\dis\f{\p}{\p t} u = (1-\theta_H) r(x)  -   d(x) - I(t) + \theta_H \;  H(x, \nabla u(x,t)) ,
\\ [5mm]
\dis \max_{0\leq x < \infty} u(x,t) =0 \qquad \forall t >0.
\eepa
\label{aah:chj2}
\eeq
As mentioned in Section \ref{sec:aah}, we find again that the fittest trait $\bar x(t)$ is characterized by
$$
0=u(\bar x(t),t), \qquad \qquad 0= r(\bar x(t)) - d(\bar x(t)) -I(t).
$$
From the definition of $I_\e(t)$, this last equality is of course equivalent to the limit $p_\e \to \delta(x-\bar x(t))$.

Again, the general theory applies here and we may conclude the proof as in Section \ref{sec:aah}: $I(t)$ is non-decreasing and thus $\bar x(t)$ is non-increasing; this proves that these two quantities have limits for $t\to \infty$, which we easily identify from the Hamilton-Jacobi  equation.
\\

With this information, notice that the dynamics of $\rho_C(t)$ can be approximated by
$$
\e \f{d  \rho_C(t) }{dt} =\rho_C(t)[ r(\bar x(t)) - d(\bar x(t)) ],
$$
and it may decrease exponentially if $r(\bar x^0) -d(\bar x^0) >0$, but eventually it will increase because we know that $r(0) -d(0) <0$.
\qed

\subsection{Asymptotic analysis: cancer cells with therapy}
\label{sec:aat}

In the case when treatment is included, the dynamics of cancer cells given by equation \eqref{eq:cancer}--\eqref{eq:rho}  is independent from healthy cells and we concentrate on them. We use the same notations as before to express the small mutations regime where we can state a clear mathematical result and write the dynamics as
\beq\begin{array}{rl}
\e \dis\f{\p}{\p t} n_C(x,t) =  &\Big [   (1-\theta_C) \; r(x)  -  d(x) -\mu_c(x) \Big]  n_C(x,t)
 +\theta_C \; \dis \int r(y) \f 1\e \wt M(y, \f{y-x}{\e}) n_C(y,t) dy.
\end{array}
\label{aa:therapy}
\eeq

\begin{theorem}[Cancer cells with therapy] We assume \eqref{as1}, \eqref{as:xr}, \eqref{aa:mutkernel}, \eqref{aa:as3} for some $0<x_0< x_c$ and that all data are Lipschitz continuous. Then, the population is monoclonal in the limit of small mutations:
$$
\f{n_C(x,t)}{\int_0^\infty n_C(x,t) dx} \underset{\e \to 0}{\longrightarrow}  \delta (x- \bar x (t)),
$$
and it selects asymptotically the fittest trait $x=x_c$
$$
\bar x(t) \underset{t \to \infty}{\longrightarrow} x_c.
$$
\label{aa:th3}
\end{theorem}

The structure of this problem is linear as in Section \ref{sec:aac}. Therefore we have to adapt the general theory of concentration as we did it for cancer cells only. There is no special difficulty compared to the proof of Theorem \ref{aa:th2} and we skip it. We just mention that, again, we cannot analyze directly the cell population density $n_c$, because it blows-up exponentially, and the natural quantity on which we can develop our method is the probability density
$$
p_\e(x,t) := \f{n_C(x,t)}{\rho_C(t)}, \qquad  \rho_C(t)  = \int_0^\infty n_C(x,t) dx.
$$

\section{A consequence for optimal therapy}
\label{sec:opt}

As we have seen in Theorem \ref{aa:th3}, the model predicts selection of a resistance gene whose level of expression $x_c$ depends on the drug dose. This leads to the question to know if cancer growth can be minimized by using an optimal dose $c$.  At this stage we only consider the case of a constant dose $c$. We analyze the particular case of coefficients used for our numerical studies in Section \ref{sec:numerics}.
\\

To state the problem, we recall that the resistant population will result in net growth (fitness) $R_c$ given by formula \eqref{as:xr}  that is
\beq
 \overline{R_c }:= r(x_c)- d(x_c) - c \mu_c(x_c) = \max_{x\geq 0} \left[r(x) -d(x)- c \mu_C(x)\right].
\label{as:fitness}
\eeq
Therapy optimization means  finding $c^*$ that minimizes $ \overline{R_c}$ over $c>0$. The most usual protocol is to maximize the dose $c$ with the only limitation of toxicity on healthy cells;  does this minimize $R_c$?

We analyze this question for the choice
\beq
R_c(x) := \frac{r_0^2}{1+x^2} - d -  \frac{\al^2}{a^2 + x^2};
\eeq
here we have set $ \al^2 := c $ to simplify the following calculations.
Natural assumptions are
\begin{align}
0 < a^2(r_0^2-d) &< \al^2,  \\
a  & < 1.
\end{align}
The first one is equivalent to  $R_c(x=0)<0$, i.e., the therapy is high enough to kill the cancer cells in their normal state $x=0$ while cancer cells by themselves have a positive growth. The second one says that the therapy acts faster on the gene expression than the natural birth process.

In order to find the maximum of $R_c(x)$, we use the variable $y=x^2$, take the derivative and set it equal to zero
$$
R'_c(y) = - \frac{r_0^2}{(1+y)^2} +  \frac{\al^2}{(a^2 + y)^2} = 0  \Longleftrightarrow \f{r_0}{1+y} = \frac{\al}{a^2 + y} .
$$
From this we obtain the condition for existence of a maximum of $R_c(y)$
\beq
a^2 r_0 < \al   \Longleftrightarrow  R'_c(0)>0,
\eeq
and the following expression for the maximum point
\beq
y_c = \frac{\al -a^2r_0  }{r_0 - \al } =
\begin{cases}
>0  \quad \text{for} \quad a^2r_0 < \al < r_0,
\\[3mm]
<0 \quad \text{otherwise}.
\end{cases}
\eeq
When $y_c >0$, we compute the maximum value of $R_c$
\beq
\overline{R_c } = \frac{(\al -r_0)^2}{1- a^2} - d.
\label{eq:rc}
\eeq

We can now conclude that:
\\ \\
$\bullet$ when $\al \geq  r_0$, then $R_c (y)$ increases from $R_c(0)<-d$ to $\overline{R_c } =-d$; this level of therapy is strong enough and resistance cannot occur.
 \\
$\bullet$  When $ 0 < a^2(r_0^2-d) \leq  \al \leq a^2 r_0 $ (this interval of $\al$ is not empty for $r_0\leq 1$ for instance), then $R_c (y)$ decreases from $\overline{R_c } =R_c(0)>-d$ to $-d$. Although weaker, this level of therapy avoids resistance to occur.
\\
$\bullet$ When $ a^2 r_0 < \al < r_0$, then $R_c (y)$ increases from $R_c(0)< 0$ to $\overline{R_c }$ and then decreasing; $\overline{R_c } $ can be positive iff $(1-a^2) r_0^2 > d$ and this can happen with $R_c(0)<0$ only when $a^2r_0 < a\sqrt{r_0^2-d}$. Resistance occurs therefore if for $\al = a\sqrt{r_0^2-d}$, then $\overline{R_c } >0$. An example is explicitly computed in Section \ref{sec:numerics}. However, since $\overline{R_c } $ in \eqref{eq:rc} is decreasing in $\al$ in the range under consideration here, the only way to avoid resistance is to increase the dose.

\section{Numerics}
\label{sec:numerics}

In this section we illustrate the analytic results obtained before by numerical simulations performed in {\sc Matlab}.
For all these computations we use 4000 points on the interval $[0,1]$. Moreover, we treat the case without mutations, i.e., $\theta_H=\theta_C=0$. The parameter $\e= 10^{-2}$ is used for the initial data which we take under the `highly concentrated form' \eqref{aa:as1}, \eqref{aa:as3} so as to obtain very localized solutions $n_{C,H}$ according to the theory in \cite{LMP}.
\\

\begin{figure}[h]
\centering
\includegraphics[width=0.45\textwidth]{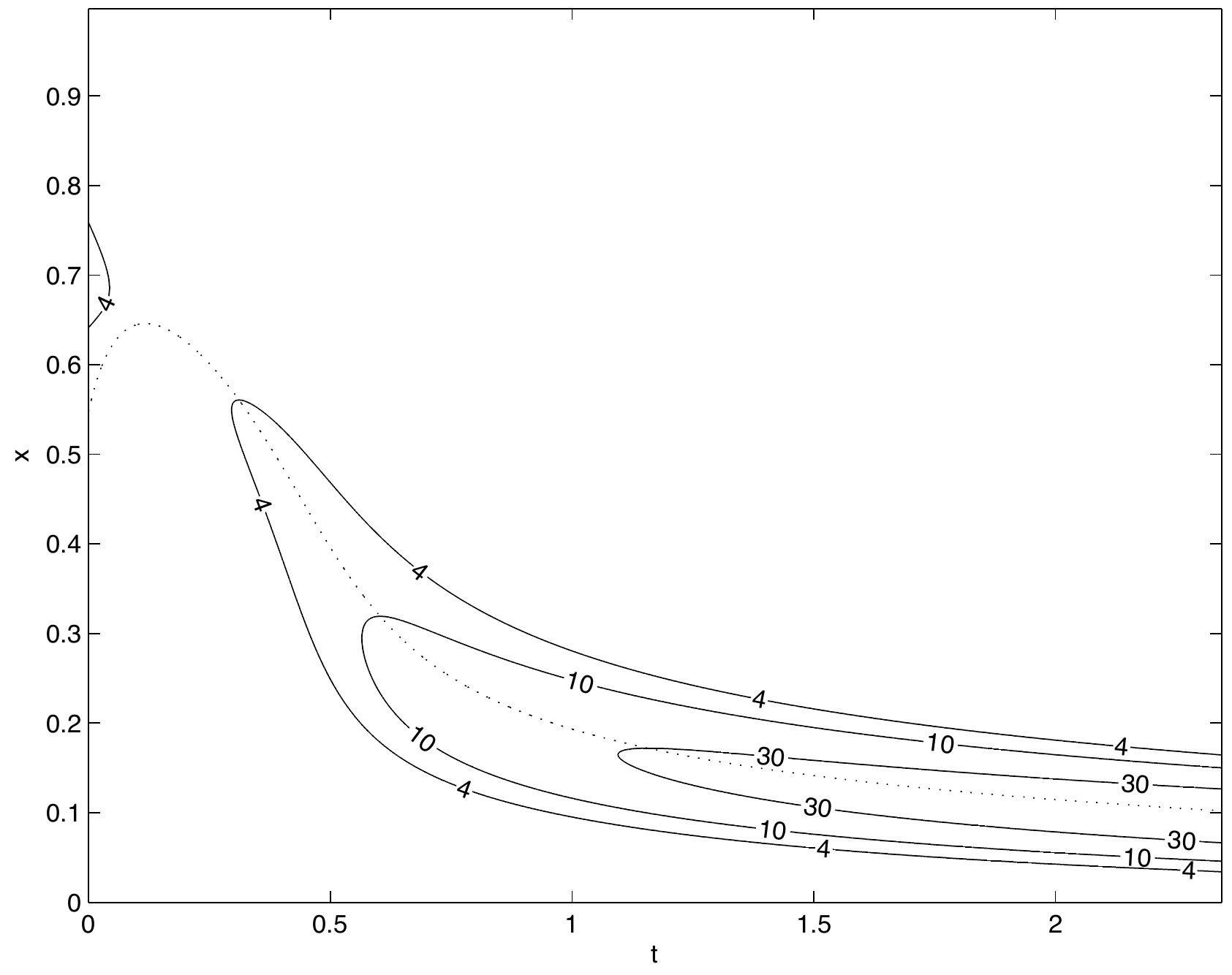}
\includegraphics[width=0.45\textwidth]{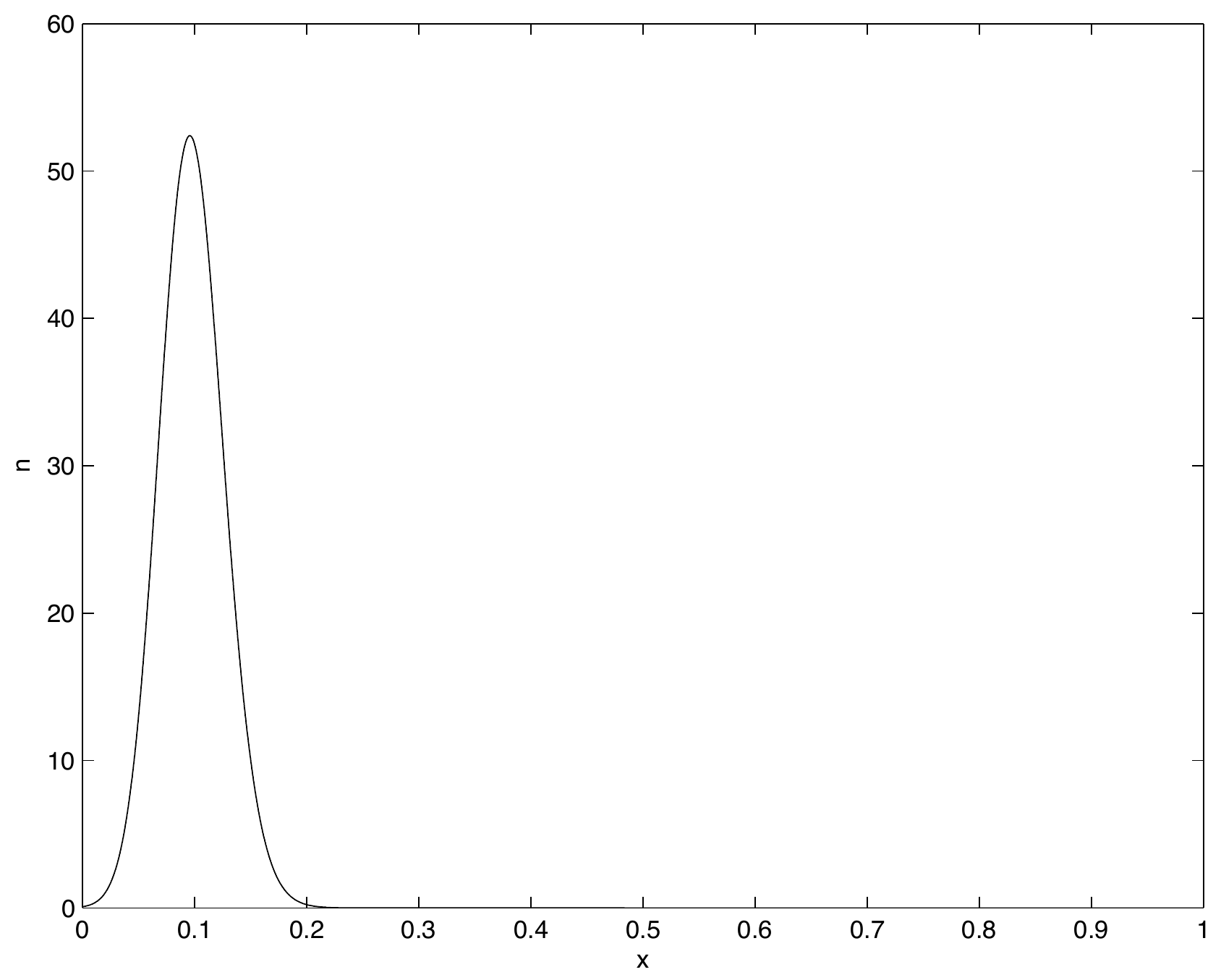}
\vspace{-4mm}
\caption{(Healthy cells) The plot to the left shows the level-sets of $n(t,x)$.
The dotted line shows the unique solution $x$ to $R(x,\bar{\rho}(t))=0$  as a function of $t$.
On the right, we display the values of $n_H(t,x)$ at the end of the computation. Notice that $n_H(t,x)$ concentrates and its maximum moves closer towards $0$ for larger time.
}
\label{fi.healthy}
\end{figure}

We first illustrate the case of healthy cells. Figure \ref{fi.healthy} shows the time dynamics of the concentration point in \re{aa:healthy} with
$$
\dis{r(x):= \f{2}{1+5x^2}}, \qquad d := 0.4
$$
and the initial data
$$
n_0(x) := C_0 \exp(-(x-0.7)^2/\e),
$$
where the constant $C_0$ is adjusted to enforce $\|n_0\|_{L^1([0,1])} = 1$. The computations were done using a implicit-explicit finite difference scheme. It can be seen that it takes very short time to adapt to $R(\bar{x}(t),\bar{\rho}(t))=0$, that $n(t,x)$ concentrates and its maximum moves towards $0$.

\begin{figure}[h]
\centering
\includegraphics[width=0.45\textwidth]{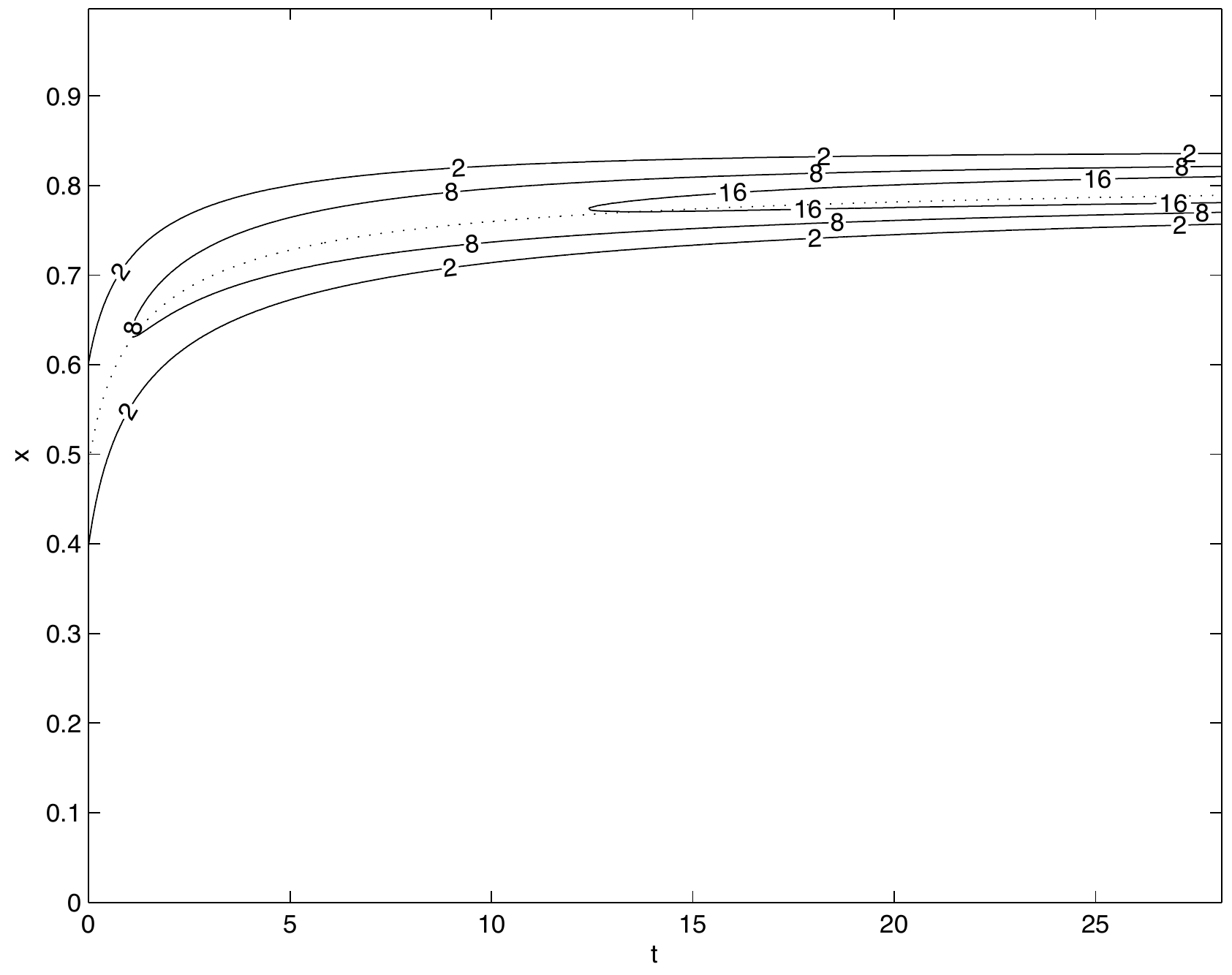}
\includegraphics[width=0.45\textwidth]{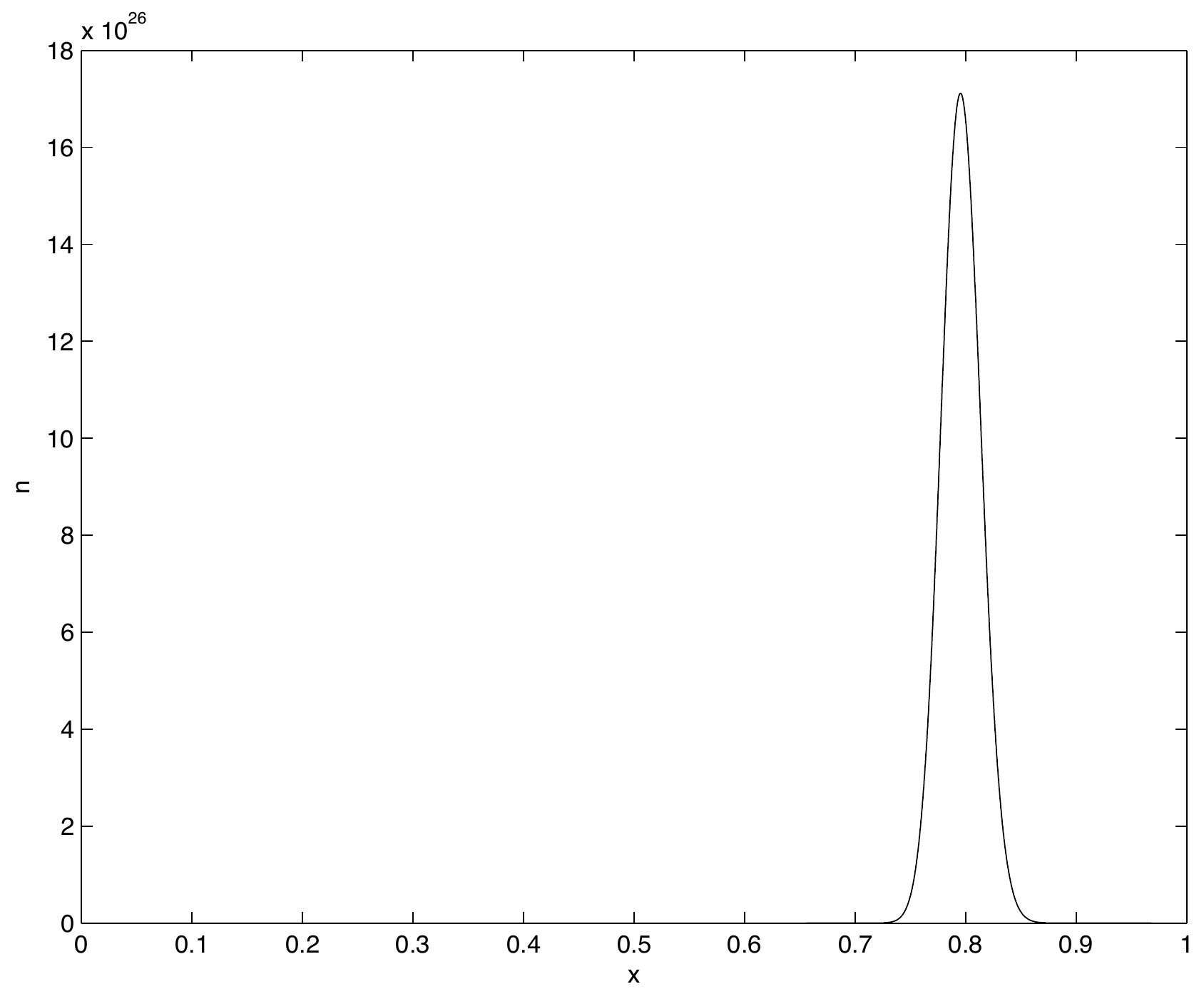}
\vspace{-4mm}
\caption{(Resistance)
The plot to the left shows the level-sets of $n_C(t,x)/\int n_C(t,x)\,dx$.
The smaller of the two positive solutions  of $r(x)-d(x) -c\mu_C(t)-I(t)=0$ is shown as  dotted line.
On the right, we plot the values of $n_C(t,x)$ at the end of the computation.
 Calculations are done for the equation \eqref{aa:therapy} on $n_C$.}
\label{fi.cancer}
\end{figure}

\begin{figure}[h!]
\centering
\includegraphics[width=0.45\textwidth]{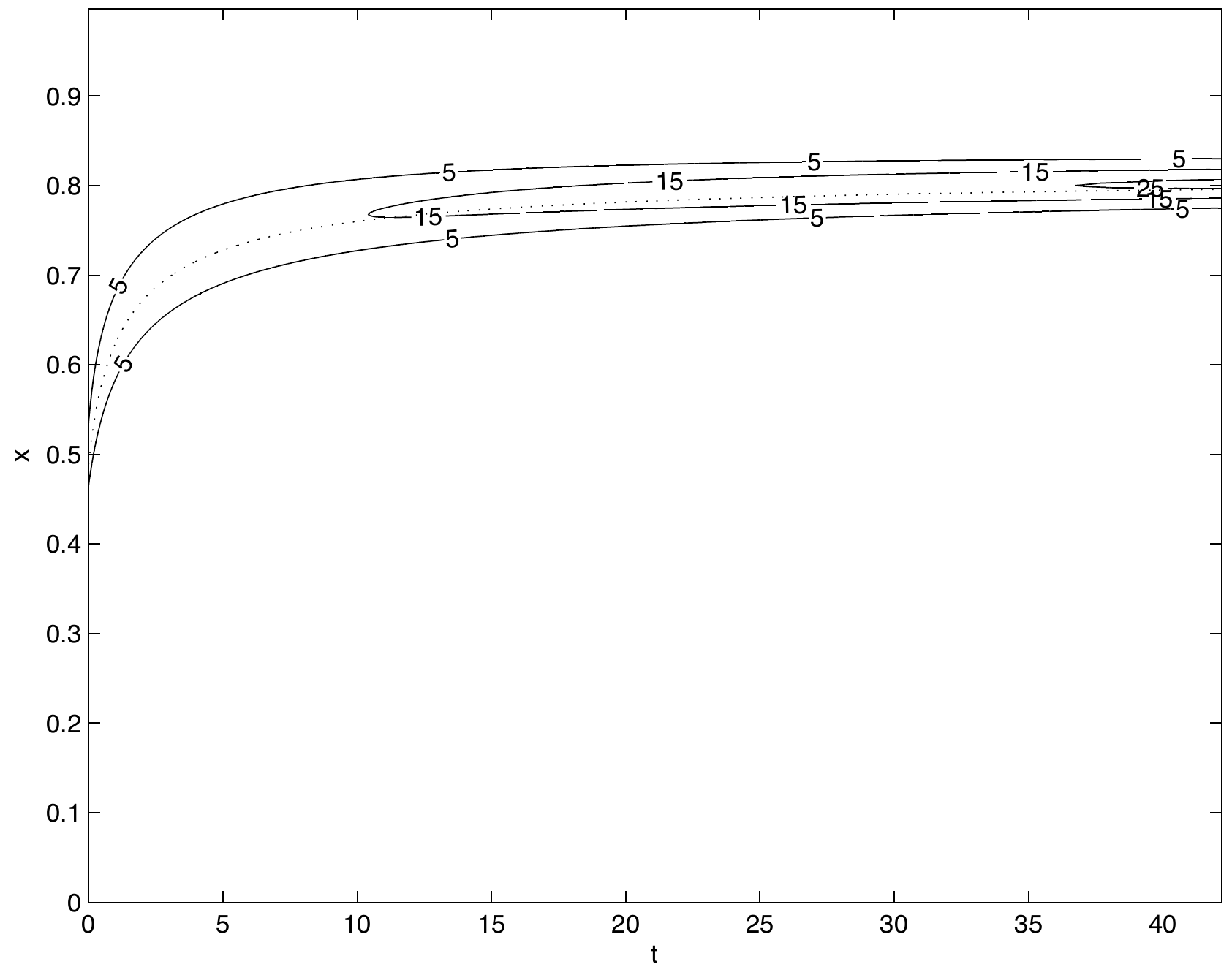}
\includegraphics[width=0.45\textwidth]{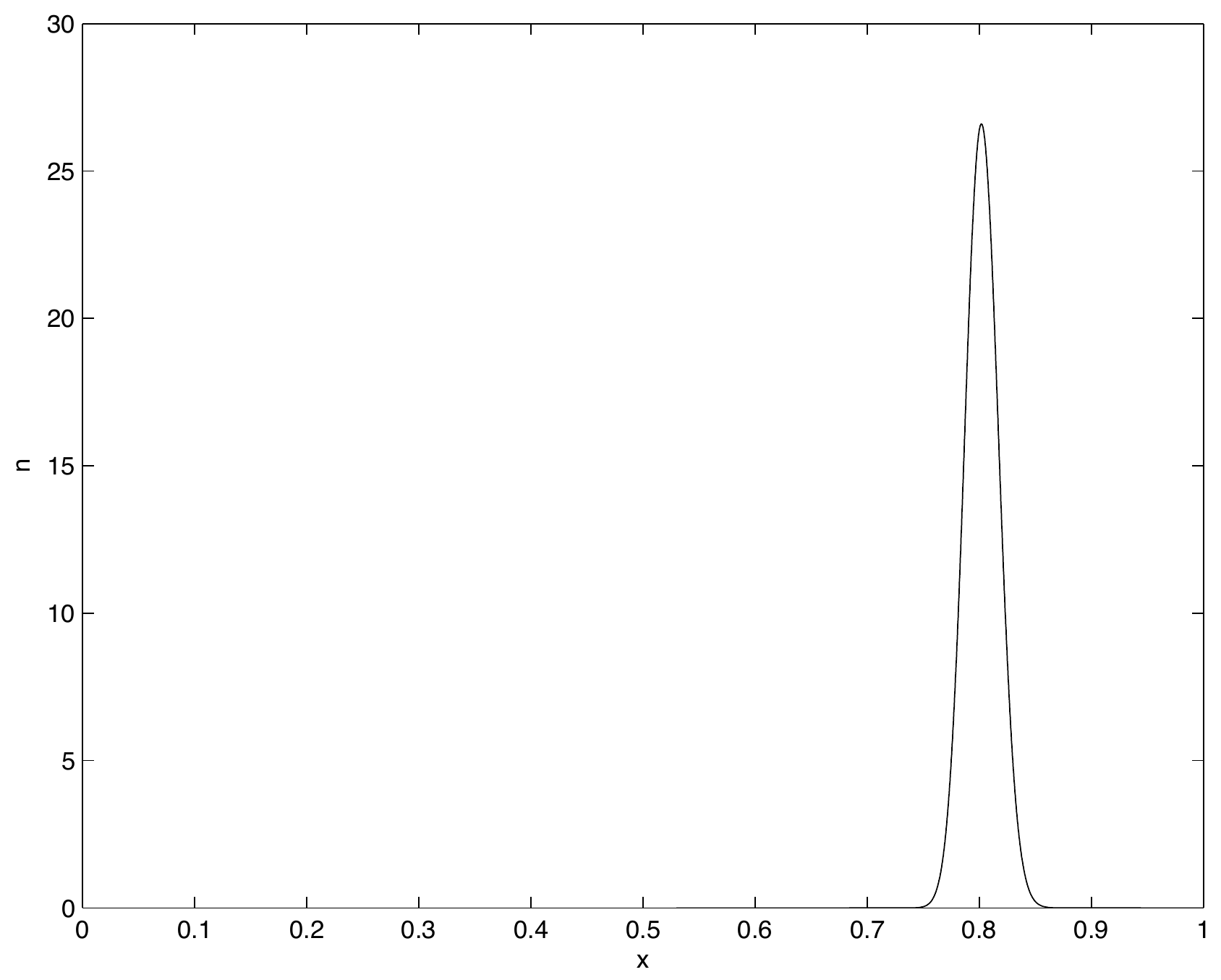}
\vspace{-4mm}
\caption{(Resistance ct'd)
The plot to the left shows the level-sets of $p$.
The smaller of the two positive solutions  of $r(x)-d(x) -c\mu_C(t)-I(t)=0$ is shown as  dotted line.
On the right, we plot the values of $p(t,x)$ at the end of the computation.
Calculations are done for the equation on $p_\e$.}
\label{fi.cancerp}
\end{figure}

Next we illustrate the case when the therapy induces resistance. Figure \ref{fi.cancer} shows the time dynamics of the concentration point in \re{aa:therapy} with
$$
\dis{r(x):= \f{1}{1+x^2}}, \qquad d := 0.245,  \qquad \dis{\mu(x):= \f{0.55^2}{0.5^2+x^2}}
$$
and initial data as before but centered at $\overline{x_0}=0.5$
$$
n_0(x) := C_0 \exp(-(x-0.5)^2/\e) .
$$
As shown in Section \ref{sec:opt}, the concentration point for large time can be computed explicitly: $$x_C = \sqrt{\f{\al-a^2}{1-\al}}=\sqrt{\f{2}{3}}\approx  0.8165.$$
Since we have exponential growth of the maximum of $n_c$, after a certain time the computations are breaking down. In order to further follow the dynamics of the maximum point $\bar{x}(t)$, we use the equation on $\dis{p_\e= \f{n_c}{\int_0^\infty n_c \,dx}}$. Results for the calculations on $p_\e$ are shown in Figure \ref{fi.cancerp}.
Figure \ref{fi.cancer} and Figure \ref{fi.cancerp} both confirm the results from theory, i.e., $r(x)-d(x) -c\mu_C(t)-I(t)=0$ (with $I(t)$ defined by \eqref{Ie}).

\section{A structured population model for healthy and tumor cells under the effects of cytotoxic and cytostatic drugs}
\label{sec:model2}

We turn to the study of the coupled dynamics of healthy and tumor cells exposed to both cytotoxic and cytostatic drugs. Our model relies on the same structured population formalism presented in Section \ref{sec:selection}. However, to allow a better qualitative agreement with biological reality, some modifications are introduced, which make the model more difficult to handle analytically. Therefore, the behavior of its solutions is studied only through numerical simulations, whose results are discussed at the end of this section.

As in Section \ref{sec:selection}, the number densities of healthy and cancer cells with resistant gene expression level $0 \leq x < \infty$ at time $t \geq 0$ are described by functions $n_H(t,x) \geq 0$ and $n_C(t,x) \geq 0$. The selection/mutation dynamics of cells under the effects of cytotoxic and cytostatic drugs are given by the system below:
\begin{eqnarray}\label{MOD2_1}
\frac{\partial}{\partial t}  n_H(x,t) & = & \overbrace{\frac{\theta_H}{1 + \alpha_H c_2(t)}\left(\int r_H(y) M(y,x) n_H(t,y) dy - r_H(x) n_H(x,t)\right)}^{\mbox{mutations and renewal}}
\nonumber\\
&+&\underbrace{\left(\frac{r_H(x)}{1 + \alpha_H c_2(t)} - d_H(x) I_H(t) \right) n_H(x,t)}_{\mbox{\mbox{growth with cytostatic therapies and death}}} - \underbrace{c_1(t) \mu_H(x) n_H(x,t)}_{\mbox{effect of cytotoxic therapies}},
\end{eqnarray}
\begin{eqnarray}\label{MOD2_2}
\frac{\partial}{\partial t} n_C(x,t) & = & \frac{\theta_C}{1 + \alpha_C c_2(t)} \left(\int r_C(y) M(y,x) n_C(t,y) dy - r_C(x) n_C(x,t)\right)
\nonumber\\
&+& \left(\frac{r_C(x)}{1 + \alpha_C c_2(t)} - d_C(x) I_C(t)\right) n_C(x,t) - c_1(t) \mu_C(x) n_C(x,t),
\end{eqnarray}
where
\begin{equation}
I_H(t) : =  a_{HH} \rho_H(t) + a_{HC} \rho_C(t), \qquad I_C(t) : = a_{CH} \rho_H(t) + a_{CC} \rho_C(t).
\end{equation}
Apart from integrals $I_H$ and $I_C$, which have a different meaning with respect to the integral $I$ considered within the previous sections, the same notations of Section \ref{sec:selection} are used in the above equations. In addition:
\\
$\bullet$ Kernel $M(y,x)$ denotes the probability that a mutation of a healthy/cancer cell with gene expression $y$ leads to a daughter cell with level $x$.
\\
$\bullet$ Functions $r_H(x)$ and $r_C(x)$ stand, respectively, for the proliferation rate of healthy and cancer cells with gene expression $x$. Factors
$$
\frac{1}{1 + \alpha_H c_2(t)} \quad \mbox{ and } \quad  \frac{1}{1 + \alpha_C c_2(t)}
$$
mimic the effects of cytostatic drugs, whose concentration at time $t$ is described by function $c_2(t)$. In fact, as previously noted, such therapies act by slowing down cell proliferation, rather than by killing cells. The average sensitivities of healthy and cancer cells to these drugs are modeled by parameters $\alpha_H, \alpha_C \in \mathbb{R}^+$.
\\
$\bullet$ Function $d_H(x)$ represents the death rate of healthy cells due both to apoptosis and deprivation of resources (e.g., oxygen and glucose) by cancer cells. Analogously, function $d_C(x)$ models the death rate of cancer cells because of the competition for space and resources with the other cells. Two main differences are here introduced with respect to the model presented in Section \ref{sec:selection}. On the one side, we assume the growth of cancer cells to be hampered by the competition for space and resources among themselves and healthy cells as well. As a result, we multiply $d_C(x)$ by $I_C(t)$. On the other side, we consider apoptosis and competition phenomena involving normal cells as mediated by interactions with the surrounding cells. Therefore, instead of multiplying function $r_H(x)$ by $\frac{1}{(1 + \varrho(t))^{\beta}}$ as in \eqref{eq:healthy}, we multiply $d_H(x)$ by $I_H(t)$. It is worth noticing that both modeling strategies mimic the same effects of net proliferation and death.
\\
$\bullet$ Function $c_1(t)$ denotes the concentration, at time $t$, of cytotoxic agents, which are assumed to kill both healthy and cancer cells with gene expression $x$ at rates $\mu_H(x)$ and $\mu_C(x)$, respectively.
\\
$\bullet$ Coefficients $a_{HC}$ and $a_{HH}$ stand, respectively, for the average interaction rates between healthy and cancer cells or among healthy cells themselves. Analogous considerations hold for $a_{CH}$ and $a_{CC}$.

As in Section \ref{sec:selection}, model \eqref{MOD2_1}, \eqref{MOD2_2} can be recast in the equivalent form given hereafter, in order to highlight the role played by the net growth rates of healthy and cancer cells, which are described by functionals $R_H(I_H,c_1,c_2,x)$ and $R_C(I_C,c_1,c_2,x)$:
\begin{eqnarray}\label{MOD3}
\frac{\partial}{\partial t}  n_H(x,t) & = & R_H(I_H(t),c_1(t),c_2(t),x) n_H(x,t) + \frac{\theta_H}{1 + \alpha_H c_2(t)} \int r_H(y) M(y,x) n_H(t,y) dy
\nonumber\\
\nonumber\\
\frac{\partial}{\partial t}  n_C(x,t) & = & R_C(I_C(t),c_1(t),c_2(t),x) n_C(x,t) + \frac{\theta_C}{1 + \alpha_C c_2(t)} \int r_C(y) M(y,x) n_C(t,y) dy,
\nonumber
\end{eqnarray}
with
\begin{eqnarray}
R_H(I_H(t),c_1(t),c_2(t),x) &:=& \frac{r_H(x)(1-\theta_M)}{1 + \alpha_H c_2(t)} - d_H(x) I_H(t) - \mu_H(x) c_1(t),
\nonumber\\
R_C(I_C(t),c_1(t),c_2(t),x) &:=& \frac{r_C(x)(1-\theta_M)}{1 + \alpha_C c_2(t)} - d_C(x) I_C(t) - \mu_C(x) c_1(t).
\nonumber
\end{eqnarray}

The following considerations and hypothesis are assumed to hold:
\\
$\bullet$ Parameters $\theta_{H,C}$ and functions $\mu_{H,C}$ satisfy the assumptions of Section \ref{sec:selection}.
\\
$\bullet$ We consider intra-population interactions as occurring at a higher rate than inter-population ones. As a consequence, the interaction rates are assumed to be real numbers such that
\begin{equation} \label{A2}
a_{HH}, a_{CC} > 0, \quad a_{HC}, a_{CH} \geq 0, \qquad a_{HH} > a_{HC}, \quad a_{CC} > a_{CH}.
\end{equation}
\\
$\bullet$ In analogy with Section \ref{sec:selection}, with the aim of translating in mathematical terms the idea that producing resistance genes implies resource allocation both for healthy and cancer cells (eventhough this is debated), we assume functions $r_H$ and $r_C$ to be decreasing
\begin{equation} \label{A3}
r_{H,C}(\cdot) > 0, \qquad r'(\cdot)_{H,C} < 0.
\end{equation}
In order to mimic the fact that mutations conferring resistance to therapies may also provide cells with stronger competitive abilities, functions $d_H$ and $d_C$ are assumed to be decreasing
\begin{equation} \label{A4}
d_{H,C}(\cdot)>0, \qquad d'(\cdot)_{H,C} < 0.
\end{equation}
\\
$\bullet$ Since we assume that cytostatic agents are designed to be more effective against cancer cells rather against healthy cells, we make the additional assumption
\begin{equation} \label{A6}
\alpha_H < \alpha_C.
\end{equation}
We analyse the model  with  these new ingredients through numerical simulations illustrating how the outputs  can be influenced by different concentrations of cytotoxic and cytostatic agents. From
a biological perspective, this means to use the present model as an \emph{in silico} laboratory to highlight some mechanisms that may play a key role in the development of cancer resistance to therapies, with the aim of providing support to the design of optimal therapeutic strategies. The following contents are organized into four subsections. Subsection \ref{subsec:Setup} presents those assumptions that define the general setup for numerical simulations. Subsection \ref{subsec:Resu1} and subsection \ref{subsec:Resu2} are devoted to study the separate effects of cytotoxic and cytostatic drugs on the dynamics of healthy and cancer cells. Finally, subsection \ref{subsec:Resu3} focuses on the combined action of these two classes of anti-cancer agents, looking for the existence of suitable doses allowing the design of optimal therapeutic strategies.
\subsection{Setup for numerical simulations}
\label{subsec:Setup}
Numerical simulations are performed  as in Section \ref{sec:numerics} with 2000 points on the interval $[0,1]$. Interval $[0,T]$ with $T=2000 dt$ is selected as time domain, where the unit time $dt$ is chosen equal to $0.1$.

We choose the initial conditions
\begin{equation}\label{condfI}
n_H(t=0,x) = n_C(t=0,x) = n^0(x) := C^0 \exp(-(x-0.5)^2/\e),
\end{equation}
where $C^0$ is a positive real constant such that
$$
\varrho_H(t=0) + \varrho_C(t=0) \approx 1.
$$
Parameter $\varepsilon$ is set equal to $0.01$ to mimic a biological scenario where most of the cells are characterized by the resistant gene expression level corresponding to $x=0.5$ at the beginning of observations.

Assumptions and definitions given hereafter are used along all simulations:
$$
M(y,x) : = C_M \exp(-(y-x)^2/(0.01)^2), \qquad C_M \int_0^1 \exp(-(y-x)^2/(0.01)^2) dx = 1, \quad \forall y \in [0,1],
$$
$$
\theta_H=\theta_C=\theta:=0.1, \qquad r_H(x) := \frac{1.5}{1+x^2}, \quad r_C(x) := \frac{3}{1+x^2}, \qquad \alpha_H:=0.01, \quad \alpha_C:=1,
$$
$$
a_{HH} = a_{CC} = a := 1, \quad a_{HC} := 0.07 \quad a_{CH} := 0.01,
$$
$$
d_H(x) := 0.5(1 - 0.1 x), \quad d_C(x) := 0.5(1 - 0.3 x),
$$
$$
\mu_H(x) := \frac{0.2}{(0.7)^2+x^2}, \quad \mu_C(x) := \frac{0.4}{(0.7)^2+x^2}.
$$
Functions $c_1(t)$ and $c_2(t)$ are assumed to be constant, i.e.
$$c_1(t):=c_1 \in \mathbb{R}^+, \qquad c_2(t):=c_2 \in \mathbb{R}^+.
$$
The values of parameters $c_1$ and $c_2$ are chosen case by case according to our aim in each subsection.

\subsection{Effects of cytotoxic drugs}
\label{subsec:Resu1}
\begin{figure}[h!]
\centerline{\includegraphics[width=1\textwidth]{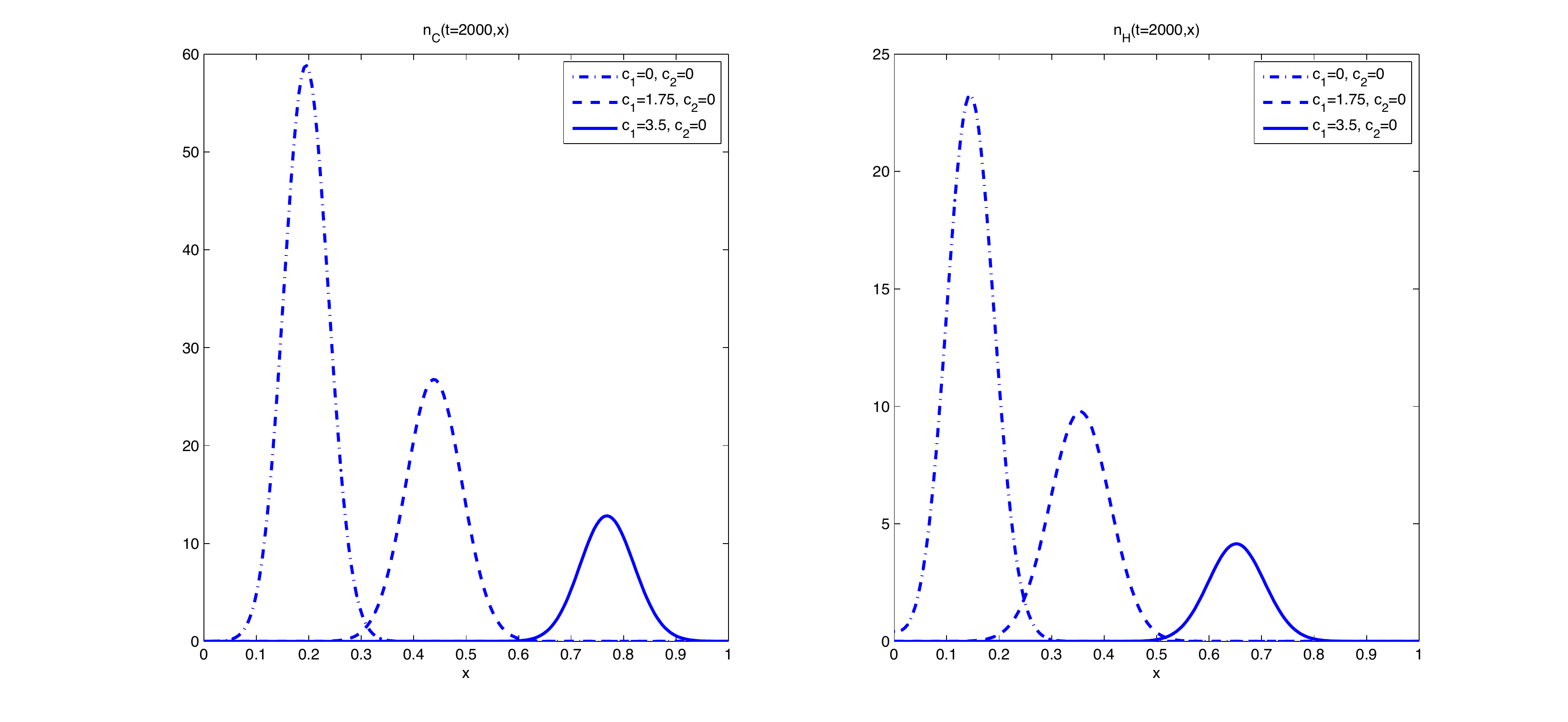}}
\vspace{-6mm}
\caption{\label{F1} (Cytotoxic drugs only) Trends of $n_C(x,t)$ (left) and $n_H(x,t)$ (right) at $t=2000$ for $c_1=0$ (dashed-dotted lines), $c_1=1.75$ (dashed lines) and $c_1=3.5$ (solid lines). In all cases, parameter $c_2$ is set equal to zero. As parameter $c_1$ increases, functions $n_C(t=2000,x)$ and $n_H(t=2000,x)$ tend to be highly concentrated around some increasing values of $x$ and their maximum values become smaller; this indicates higher resistance with higher doses of drug.}
\end{figure}

At first, we study the effects of cytotoxic agents only on the dynamics of healthy and cancer cells. With this aim, we run simulations for different values of parameter $c_1$ with $c_2=0$. The obtained results are summarized by Figure \ref{F1}, which depicts a scenario that is in good agreement with clinical observations.

In absence of therapeutic agents (i.e., when $c_1=c_2=0$), we observe, at the end of the simulations, the selection for those cells that are characterized by a low expression level of the resistant genes and thus, due to assumption \eqref{A3}, by a strong proliferative potential. On the other hand, as long as $c_1$ increases, the number of cancer cells inside the sample at the end of the observation window becomes smaller but, due to the typical side-effects of cytotoxic agents, even the number of healthy cells decreases. Furthermore, Figure \ref{F1} highlights how cytotoxic agents alone favor the selection of resistant clones. In fact, the average resistant gene expression level of cells tends to increase with the concentrations of therapeutic agents.

\subsection{Effects of cytostatic drugs}
\label{subsec:Resu2}
\begin{figure}[h!]
\centerline{\includegraphics[width=1\textwidth]{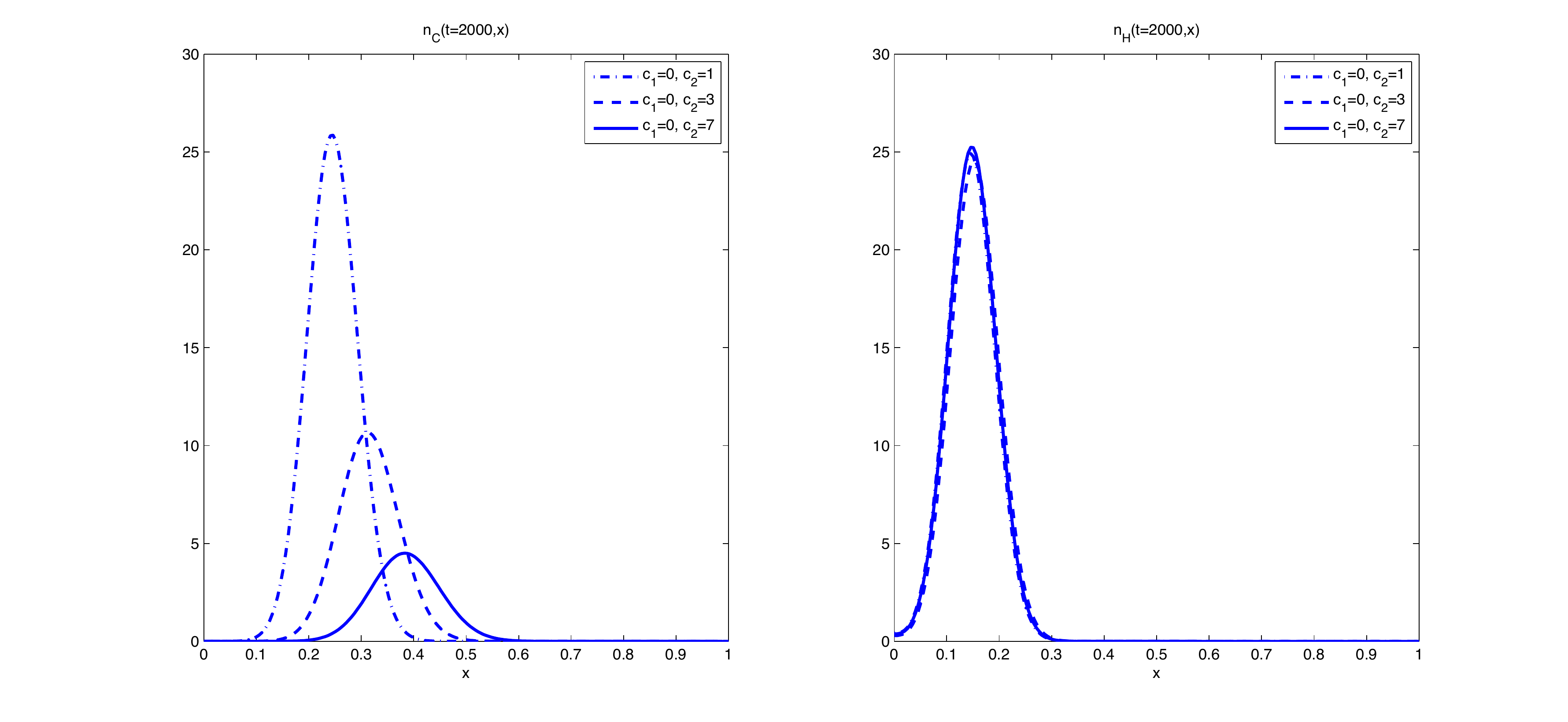}}
\vspace{-6mm}
\caption{\label{F2}(Cytostatic drugs only) Trends of $n_C(x,t)$ (left) and $n_H(x,t)$ (right) at $t=2000$ for $c_2=1$ (dashed-dotted lines), $c_2=3$ (dashed lines) and $c_2=7$ (solid lines). In all cases, parameter $c_1$ is set equal to zero. Increasing values of parameter $c_2$ lead the qualitative behavior of function $n_C(t=2000,x)$ to become closer to the one of $n_C(t=0,x)$. On the other hand, the trend of function $n_H(t=2000,x)$ remains basically unaltered.}
\end{figure}

Focusing on the action of cytostatic agents alone, we run simulations increasing parameter $c_2$ and keeping $c_1=0$. The results presented in Figure \ref{F2} highlight the capability of the present model to mimic the effects typically induced by cytostatic agents on cells dynamics. In fact, under the considered values of $c_2$, the dynamics of healthy cells is kept unaltered while the proliferation of cancer cells is reduced as long as the drug concentration increases. Furthermore, it should be noted that increasing values of $c_2$ tend to slow down the dynamics of cancer cells, i.e., for larger values of $c_2$, function $n_C(x,t)$ at the end of simulations stays closer to the initial data $n_C(t=0,x)$.

\subsection{Combined effects of cytotoxic and cytostatic drugs}
\label{subsec:Resu3}
\begin{figure}[h!]
\centerline{\includegraphics[width=1\textwidth]{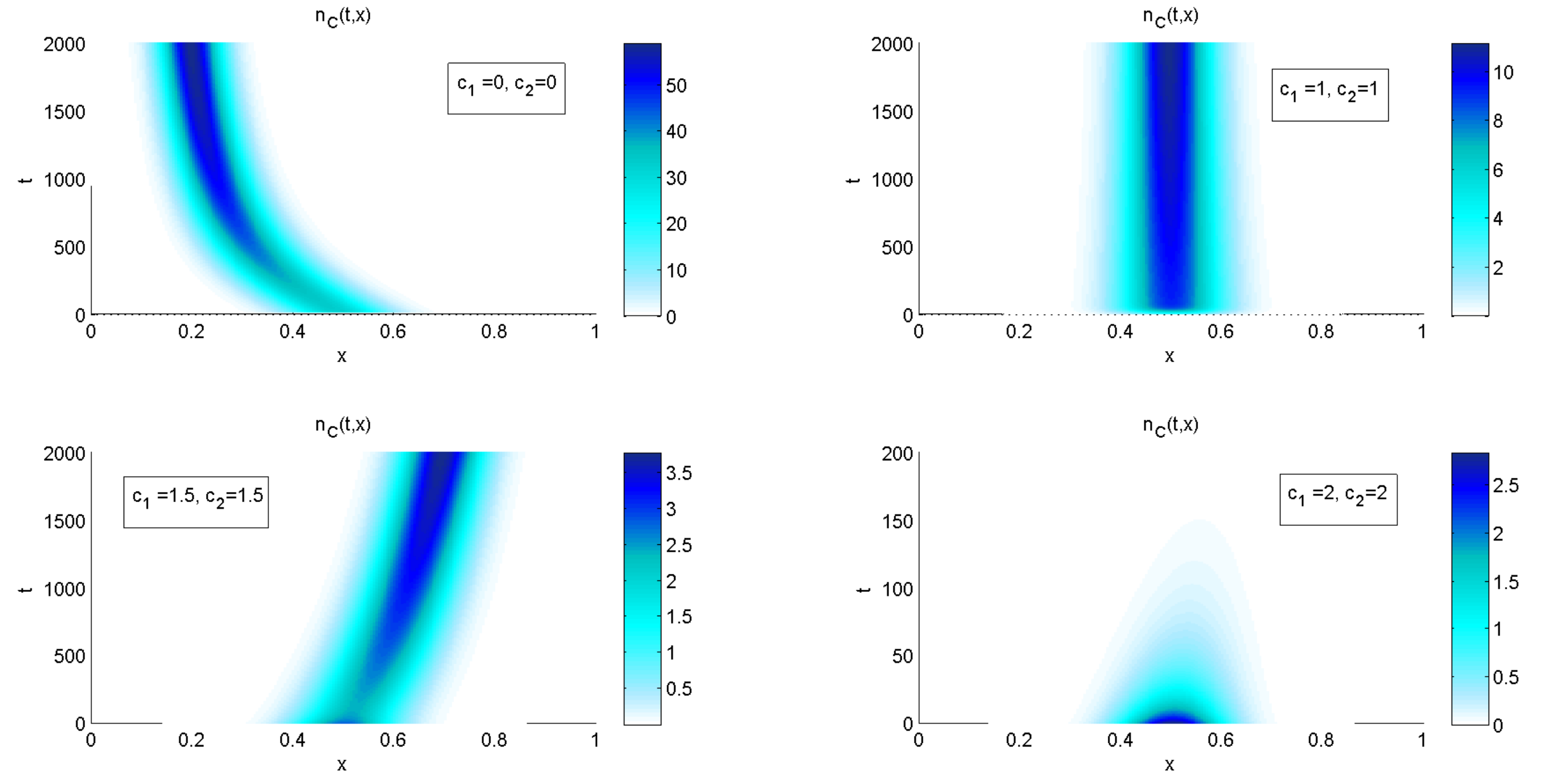}}
\vspace{-6mm}
\caption{(\label{F3}Cytotoxic and cytostatic drugs) Dynamics of $n_C(x,t)$ for $c_1=c_2=0$ (top-left), $c_1=c_2=1$ (top-right), $c_1=c_2=1.5$ (bottom-left) and $c_1=c_2=2$ (bottom-right). As long as parameters $c_1$ and $c_2$ increase, the maximum value of $n_C(t=2000,x)$ becomes smaller so that, under the choice $c_1=c_2=2$, function $n_C(x,t)$ tends to zero across time.}
\end{figure}

\begin{figure}[h!]
\centerline{\includegraphics[width=1\textwidth]{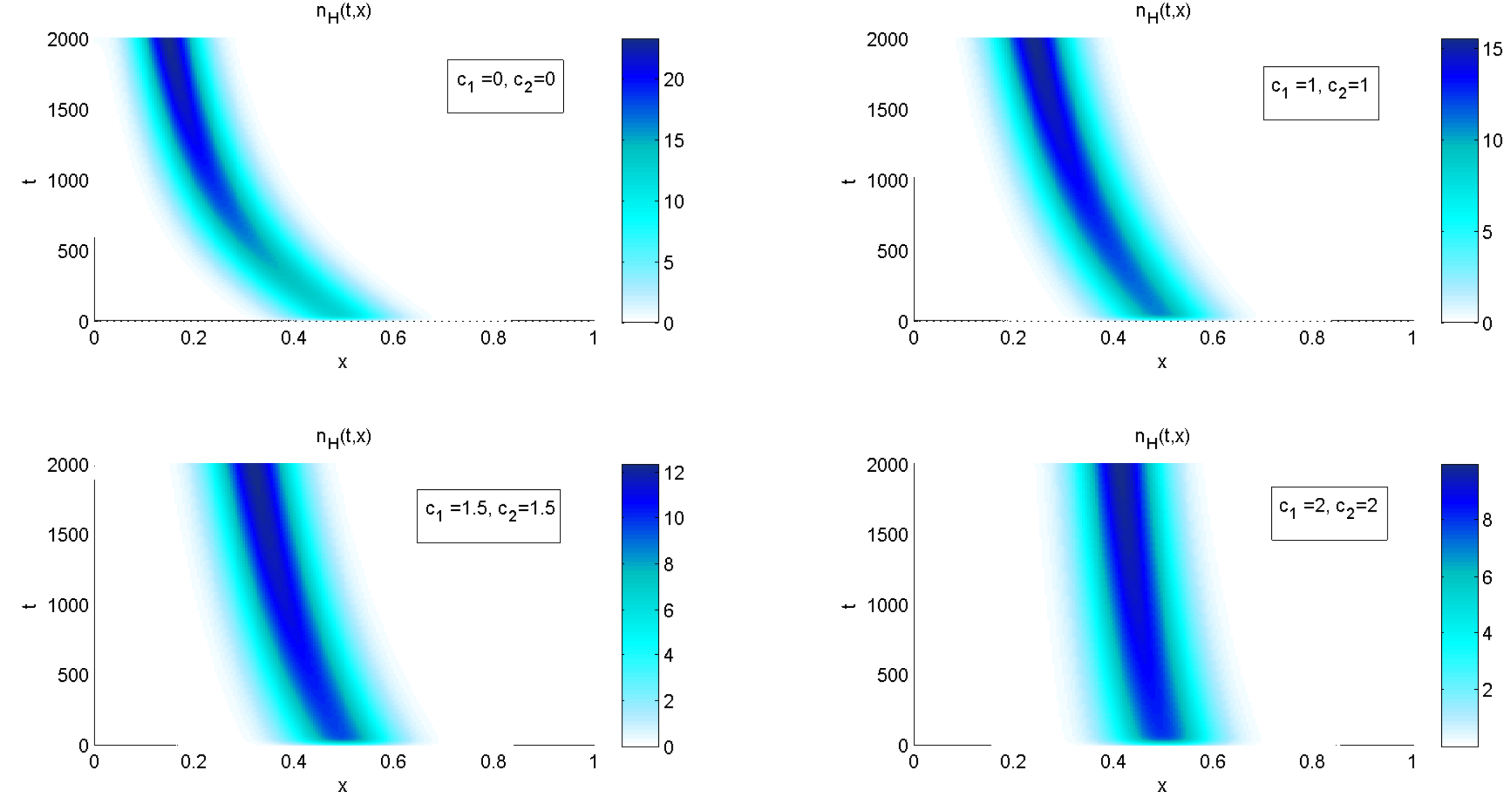}}
\vspace{-6mm}
\caption{\label{F4}(Cytotoxic and cytostatic drugs) Dynamics of $n_H(x,t)$ for $c_1=c_2=0$ (top-left), $c_1=c_2=1$ (top-right), $c_1=c_2=1.5$ (bottom-left) and $c_1=c_2=2$ (bottom-right). As long as parameters $c_1$ and $c_2$ increase, the maximum value of $n_H(t=2000,x)$ becomes smaller but about one half of the healthy cells is still alive at the end of computations.}
\end{figure}

Finally, we are interested in analyzing the effects on cell dynamics caused by the simultaneous action of cytotoxic and cytostatic drugs. As a result, we perform simulations with various values of  both $c_1$ and $c_2$. The obtained results are summarized by Figure \ref{F3} and Figure \ref{F4}. They show that there are values of $c_1$ and $c_2$ leading to extinction of the cancer cells while keeping alive about one half of the healthy cells at the end of computations. This is consistent with experimental observations suggesting that different therapies in combination can avoid the emergence of resistance and minimize side-effects on healthy cells. Furthermore, it should be noted that such values of $c_1$ and $c_2$ are smaller than the highest ones considered in previous simulations, which do not even allow a complete eradication of cancer cells  (see Figure \ref{F1} and Figure \ref{F2}). Therefore, these results also suggest that optimized anti-cancer treatments can be designed making use of proper combinations between cytotoxic and cytostatic agents, thus supporting the idea that looking for protocols based on treatment optimization can be a more effective strategy for fighting cancer rather than only using high drug doses.

\section{Conclusion and perspectives}
\label{sec:perspectives}

Based on a model adapted from ecology, we have presented a mathematical analysis of resistance mechanisms to drugs under the assumption that resistance is induced by adaptation to drug environmental pressures. We have used physiologically structured equations where the variable represents the gene expression levels. Theories developed in other contexts of population biology and Darwinian evolution underlay our analysis and finding that therapy can trigger emergence of a resistance gene and result in disease relapse.
\\

Importantly, results suggest the feasibility of optimized anti-cancer treatments, based on specific combinations of cytotoxic and cytostatic drugs (see Section \ref{sec:model2}). On the other hand, if only cytotoxic drugs are available, then (not surprisingly) the optimal therapy is always to administer the maximal tolerated dose (see Section \ref{sec:opt}). This scenario is somewhat different from that studied in \cite{GSGF}, where the resistant and susceptible populations are in competition, a polymorphism that cannot occur in our model because it contains a single environmental variable (Gause competitive exclusion principle). It is worth noting that we have analyzed the possibility for optimization only in the restrictive case of a constant level of therapy and a possible future extension would be to modulate dose with time \cite{FooM2}.
\\

The model is very simple at this stage and several directions for improvements are possible. As mentioned above, because there is only a single environmental unknown (the drug pressure and only one trait depending on it), polymorphism is not possible and thus a more realistic description of the micro-environment of the cell populations is needed (see \cite{jc} and the references therein). 
Another research direction would be to include several qualitatively different drug resistance mechanisms, together with realistic, targeted, and clinically acceptable multi-drug therapies.
\\

\bigskip

\noindent{\em Acknowledgments.}  The authors wish to thank J.-P. Marie (H\^opital Saint-Antoine, Paris) for many helpful discussions and insights on the modeling issues. JC, MEH and BP are supported by a PEPii grant from the Centre National de la Recherche Scientifique.

%
%
%

\appendix
\section{Matlab code for Figure \ref{fi.healthy}}
\begin{verbatim}
% solves  equation
% \partial_t n = [2./((1.+5.*x.^2).*(1 + rho))-0.4]n
%------------------------------------------------------------
%------------------------------------------------------------
clear all;
close all;
m=4000;         %number of points for space grid

eps=.01;        %epsilon
rhoin=1.;       %initial rho

dx=1/m; x = (0:m)*dx;  %space grid (0,1/m,2/m,...,1)
dt=25.*dx^2/eps;       %time step

dtsave = dt;            %initialize array to save time steps

tfinal = 15000 * dt;    %final time


weights = ones(numel(x));  %weights for trapez rule integration
weights(1) = 0.5;
weights(end) = 0.5;

nd=exp(-(x-.7).^2/eps);    %initial data
nd=rhoin*nd/(dx*sum(nd));  %initial data normalized to have mass rhoin
ndd=nd;R=x;                %initialize ndd and R
count = 0;

temps=0;
savend(1,:)=x(1:10:m);      %initialize arrary for saving nd at every 10th point in space

while temps < tfinal
    count = count +1
    dt = min(dt, tfinal -temps);
    rho = dx * trapz(nd)
    savend(count,:)=nd(1:10:m);
    arrayt(count) = temps;
    arrayR(count) = ((4-rho)/(5+5*rho))^(1/2);

    R=2./((1.+5.*x.^2).*(1 + rho))-0.4;

    Rp=max(0,R)/eps; Rm=min(0,R)/eps;       %positive and negative part of R

    ndd = ((1+dt*Rp).*nd )./(1-dt*Rm);      %calculate n(t+dt) from n(t)

    nd=ndd;
    temps = temps + dt;

end

figure
[X,Y] = meshgrid(arrayt,x(1:10:m));

[C,h]=contour(X,Y,savend',[4 10 30],'LineColor','black');
text_handle = clabel(C,h);
hold on;
plot(arrayt,arrayR,':k');
xlabel('t')
ylabel('x')

figure;
plot(x,nd,'k');
xlabel('x')
ylabel('n')

\end{verbatim}

\section{Matlab code for Figure \ref{fi.cancer} }
\begin{verbatim}
% solves  equation
% \partial_t n = [1./(1.+x.^2)-dea-(al^2)./(aa^2+x.^2)]n
%-------------------------------------------------------------
% ------------------------------------------------------------
clear all;
close all;
m=4000;         %number of points for space grid

eps=.01;        %epsilon
rhoin=1.;       %initial rho

dx=1/m; x = (0:m)*dx;   %space grid (0,1/m,2/m,...,1)

dt=4500*dx^2/eps;       %time step

dtsave = dt;            %initialize array to save time steps

tfinal = 1000 * dt;     %final time

weights = ones(numel(x));   %weights for trapez rule integration
weights(1) = 0.5;
weights(end) = 0.5;

nd=exp(-(x-.5).^2/eps);     %initial data
nd=rhoin*nd/(dx*sum(nd));   %initial data normalized to have mass rhoin

ndd=nd;R=x;                 %initialize ndd and R
nd0=nd;
count = 0;
temps=0;
savendrenorm(1,:)=x(1:10:m);    %initialize arrary for saving nd at every 10th point in space
% saveR(1,:)=x(1:10:m);           %initialize arrary for saving R at every 10th point in space
al=.55
aa=.5
dea=.245

R=1./(1.+x.^2)-dea-(al^2)./(aa^2+x.^2);

while temps < tfinal
    count = count +1

    nd = exp(R*temps/eps).*nd0; %exponential formula to solve the ODE

    pp = nd/(dx * trapz(nd));
    Ik= dx * trapz(pp.*R);

    dt = min(dt, tfinal -temps);
    rho = dx * trapz(nd)

    savendrenorm(count,:)=nd(1:10:m)/(dx * trapz(nd));
    arrayt(count) = temps;
        %calculating explicitly the position of the concentration point
        aq=(Ik+dea)
        bq=-1+al^2+(Ik+dea)*(1+aa^2)
        cq=-aa^2+al^2+(Ik+dea)*aa^2

        xq2=(-bq-(bq^2-4*aq*cq)^(1/2))/(2*aq)
        xq = (xq2)^(1/2)
        arrayR(count) = xq;

    temps = temps + dt;


end

figure
[X,Y] = meshgrid(arrayt,x(1:10:m));

[C,h]=contour(X,Y,savendrenorm',[2 8 16],'LineColor','black');
text_handle = clabel(C,h);
hold on;
plot(arrayt,arrayR,':k');
xlabel('t')
ylabel('x')

figure;
plot(x,nd,'k');
xlabel('x')
ylabel('n')
\end{verbatim}

\end{document}